\documentclass[12pt]{amsart}
\usepackage{amsmath,amscd,amssymb,amsfonts}
\newcommand{\vers}{Jan. 3, 2004, v.2.1}
\newcommand{\scirc}{\raise.2ex\hbox{${\scriptstyle\circ}$}}
\newcommand{\mtims}{\hbox{$\times$}}
\newcommand{\msum}{\hbox{$\sum$}}

\newcommand{\mcup}{\hbox{$\bigcup$}}
\newcommand{\mcap}{\hbox{$\bigcap$}}
\newcommand{\bA}{{\mathbb A}}
\newcommand{\bC}{{\mathbb C}}
\newcommand{\bD}{{\mathbb D}}
\newcommand{\bP}{{\mathbb P}}
\newcommand{\bQ}{{\mathbb Q}}
\newcommand{\bZ}{{\mathbb Z}}
\newcommand{\boP}{{\mathbf P}}
\newcommand{\boR}{{\mathbf R}}
\newcommand{\cA}{{\mathcal A}}
\newcommand{\cC}{{\mathcal C}}
\newcommand{\cE}{{\mathcal E}}
\newcommand{\cF}{{\mathcal F}}
\newcommand{\cG}{{\mathcal G}}
\newcommand{\cH}{{\mathcal H}}
\newcommand{\cI}{{\mathcal I}}
\newcommand{\cL}{{\mathcal L}}
\newcommand{\cO}{{\mathcal O}}
\newcommand{\cP}{{\mathcal P}}
\newcommand{\cV}{{\mathcal V}}
\newcommand{\cX}{{\mathcal X}}
\newcommand{\cZ}{{\mathcal Z}}

\newcommand{\oj}{\bar{j}}
\newcommand{\oS}{\bar{S}}

\newcommand{\oY}{\bar{Y}}
\newcommand{\oZ}{\bar{Z}}
\newcommand{\tf}{\widetilde{f}}
\newcommand{\tH}{\widetilde{H}}
\newcommand{\tL}{\widetilde{L}}
\newcommand{\tV}{\widetilde{V}}

\newcommand{\tY}{\widetilde{Y}}
\newcommand{\tcX}{\widetilde{\mathcal X}}
\newcommand{\Gr}{\text{\rm Gr}}
\newcommand{\prim}{\text{\rm prim}}
\newcommand{\pr}{\text{\rm pr}}
\newcommand{\BM}{\text{\rm BM}}
\renewcommand{\Im}{\text{\rm Im}}
\newcommand{\Ker}{\text{\rm Ker}}
\newcommand{\codim}{\text{\rm codim}}
\newcommand{\rank}{\text{\rm rank}}
\newcommand{\Sing}{\text{\rm Sing}}
\newcommand{\Coker}{\text{\rm Coker}}
\newcommand{\Ext}{\text{\rm Ext}}
\newcommand{\Hom}{\text{\rm Hom}}
\newcommand{\cHom}{{{\mathcal H}om}}
\newcommand{\van}{\text{\rm van}}
\newcommand{\reg}{\text{\rm reg}}
\newcommand{\can}{\text{\rm can}}
\newcommand{\var}{\text{\rm var}}
\newcommand{\Tr}{\text{\rm Tr}}

\begin{document}
\title[Monodromy of a family of hypersurfaces]
{Monodromy of a family of hypersurfaces\\
containing a given subvariety}
\author{Ania Otwinowska}
\address{Laboratoire de Math\'ematiques\\
Universit\'e Paris-Sud - B\^at 425\\
91405 Orsay Cedex, France}
\email{Ania.Otwinowska@math.u-psud.fr}
\author{Morihiko Saito}
\address{RIMS Kyoto University\\
Kyoto 606--8502\\
Japan}
\email{msaito@kurims.kyoto-u.ac.jp}
\date{\vers}
\begin{abstract}
For a subvariety of a smooth projective variety, consider the
family of smooth hypersurfaces of sufficiently large degree
containing it, and take the quotient of the middle cohomology
of the hypersurfaces by the cohomology of the ambient variety
and also by the cycle classes of the irreducible components
of the subvariety.
Using Deligne's semisimplicity theorem together with
Steenbrink's theory for semistable degenerations,
we give a simpler proof of the first author's theorem
(with a better bound of the degree of hypersurfaces)
that this monodromy representation is irreducible.
\end{abstract}
\maketitle

\centerline{\bf Introduction}

\bigskip\noindent
It is well known after Noether and Lefschetz that, for a
{\it general} smooth hypersurface
$ X $ in
$ \bP_{\bC}^{3} $, the subspace of Hodge cycles in the
middle cohomology
$ H^{2}(X,\bQ) $ is generated by that of
$ \bP_{\bC}^{3} $ (i.e. it is
$ 1 $-dimensional) if the degree of the hypersurface is at least
$ 4 $.
This follows from the irreducibility of the monodromy
representation on the primitive cohomology of the hypersurfaces.
Sometimes we want to consider a family of hypersurfaces
containing a given closed subvariety
$ Z $, and ask if an analogue of the above assertion holds.

More generally, let
$ Y $ be an irreducible smooth complex projective
variety embedded in a projective space, and
$ Z $ be a (possibly reducible) closed subvariety of
$ Y $.
Let
$ Z_{\{j\}} = \{x \in Z : \dim T_{x}Z = j\} $, where
$ T_{x}Z $ denotes the Zariski tangent space.
In this paper we assume

\medskip\noindent
(0.1)\quad
$ \dim Z_{\{j\}} + j \le \dim Y - 1 \,\, $ for any
$ j \le \dim Y $.

\medskip\noindent
In particular,
$ 2\dim Z \le \dim Y - 1 $.
Let
$ \cI_{Z} $ be the ideal sheaf of
$ Z $, and
$ \delta $ be a positive integer such that
$ \cI_{Z}(\delta) \,(:= \cI_{Z}\otimes \cO_{Y}(\delta)) $ is
generated by global sections, where
$ \cO_{Y}(i) $ denotes the restriction of
$ \cO(i) $ on the ambient projective space.
S. Kleiman and A. Altman [13] then proved that condition (0.1)
implies the existence of a smooth hypersurface section
$ X $ of degree
$ d $ of
$ Y $ containing
$ Z $ for any
$ d > \delta $ (where their condition for
$ \delta $ is slightly different from ours).
Actually, condition (0.1) is also a necessary condition,
see [19] (or (1.4) below).

Let
$ m = \dim X $, and
$ H^{m}(X,\bQ)^{\van} $ denote the orthogonal complement of
the injective image of
$ H^{m}(Y,\bQ) $ in
$ H^{m}(X,\bQ) $, which is called the vanishing part
or the vanishing cohomology,
because it is generated by the vanishing cycles of a Lefschetz
pencil using the Picard-Lefschetz formula.
Let
$ H^{m}(X,\bQ)_{Z}^{\van} $ denote the subspace of
$ H^{m}(X,\bQ)^{\van} $ generated by the cycle classes of
the maximal dimensional irreducible components of
$ Z $ modulo the image of
$ H^{m}(Y,\bQ) $ (using the orthogonal decomposition)
if
$ m = 2\dim Z $, and
$ H^{m}(X,\bQ)_{Z}^{\van} = 0 $ otherwise.
Let
$ H^{m}(X,\bQ)_{\bot Z}^{\van} $ be its orthogonal complement
in
$ H^{m}(X,\bQ)^{\van} $.
The first author [19] made the following

\medskip\noindent
{\bf 0.2.~Conjecture.}
Assume
$ \deg X \ge \delta + 1 $.
Then the monodromy representation on
$ H^{m}(X,\bQ)_{\bot Z}^{\van} $ for the family of hypersurface
sections
$ X $ containing
$ Z $ as above is irreducible.

\medskip
In the case
$ Z $ is smooth, an easy proof was given by C. Voisin,
see (1.6).
In general, using a degeneration argument inspired by
[10], [16], the first author proved the following.

\medskip\noindent
{\bf 0.3.~Theorem} [19].
{\it There exists a positive real number
$ C $ such that Conjecture {\rm (0.2)} holds if
$ \deg X \ge C(\delta+1) $.
}

\medskip
This was used in an essential way for the proof of
the main theorem in [20], which implies the Hodge conjecture
for a hypersurface section of sufficiently large degree belonging
to some open subset of an irreducible component of the
Noether-Lefschetz locus of low codimension, and whose argument
generalizes that of [18].
In the proof of Theorem (0.3), however, an asymptotic argument
is used, and
$ C $ can be quite large.
In this paper, we give a simpler proof of the theorem below
by using the theory of nearby cycles ([5], [24])
together with Hodge theory ([2], [4], [22], [25], [29])
and some cohomological properties of Lefschetz pencils
([6], [12], [14], [15], [28]).

\medskip\noindent
{\bf 0.4.~Theorem.}
{\it Let
$ \delta $ be as above, and
$ d $ be a positive integer.
Assume either
$ d \ge 2+\delta $
or a general hypersurface section of degree
$ d-\delta $ of
$ Y $ has a nontrivial differential form of the highest
degree.
Then Conjecture {\rm (0.2)} holds with the assumption
replaced by
$ \deg X = d $.
}

\medskip
The first assumption
$ d \ge \delta+2 $ may be replaced with a weaker condition that
the
$ 2 $-jets at each point are generated by the global sections
of
$ \cO_{Y}(d-\delta) $.
The second assumption on a differential form is stable by
hypersurface sections, and it is satisfied for any
$ d \ge \delta+1 $ if
$ Y $ has a nontrivial differential form of the highest
degree.
The proof of Theorem (0.4) uses Deligne's semisimplicity theorem
[4] which implies that irreducibility is equivalent to
indecomposability.
We take a semistable degeneration as in [10], and calculate
the graded pieces of the weight filtration on the nearby cycles
using [24].
Then we can proceed by induction on
$ m $ showing the nontriviality of certain extension classes
by calculating the cohomology of a Lefschetz pencil or a
non Lefschetz fibration.
In the non Lefschetz case, we use a special kind of
degeneration generalizing a construction in the surface
case in [19] (see (2.7) below) if the condition
$ d \ge \delta+2 $ is satisfied, and use Hodge theory (see (3.2),
(3.3) below) if the second assumption on a differential form
is satisfied.
In the Lefschetz case, the argument is rather easy,
see (2.3), (2.5) below.
(We can generalize Theorem (0.4) by replacing
$ \cO_{Y}(\delta) $,
$ \cO_{Y}(d-\delta) $ with two ample line bundles
$ \cL_{1} $,
$ \cL_{2} $ satisfying appropriate conditions even if they are
linearly independent in the Picard group of
$ Y $ tensored with
$ \bQ $, see (4.8) below.)

By a standard argument, Theorem (0.4) implies

\medskip\noindent
{\bf 0.5.~Corollary.}
{\it Under the assumption of Theorem {\rm (0.4)}, assume
further that
$ m = 2\dim Z $, and the vanishing cohomology of a general
hypersurface section of degree
$ d $ does not have a Hodge structure of type
$ (m/2,m/2) $.
Then the Hodge cycles
in the middle cohomology of a general hypersurface section
$ X $ of degree
$ d $ of
$ Y $ containing
$ Z $ are generated by the image of the Hodge cycles on
$ Y $ together with the cycle classes of the irreducible
components of
$ Z $.
In particular, the Hodge conjecture for
$ X $ is reduced to that for
$ Y $.
}

\medskip
In Sect.~1, we review the theory of hypersurface sections
containing a subvariety (this is mainly a reproduction from
[19]).
In Sect.~2, we prove the nonvanishing of some extension
classes using a topological method, and
in Sect.~3, we do it using a Hodge-theoretic method.
In Sect.~4, we study the graded pieces of the weight
filtration on the nearby cycles, and prove Theorem (0.4).

\medskip
In this paper, a variety means a (not necessarily reduced nor
irreducible) separated scheme of finite type over
$ \bC $, and a point of a variety means a closed point.
We say that a member of a family parametrized by the
points of a variety is general if the point corresponding to
the member belongs to some dense open subvariety of the variety.

\bigskip\bigskip
\centerline{\bf 1. Hypersurfaces Containing a Subvariety}

\bigskip\noindent
In this section, we review some basic facts in the theory of
hypersurface sections containing a subvariety, which will
be needed in the proof of Theorem (0.4).
We mainly reproduce the arguments in [19], Section 1,
see also [13].
Our argument works only in the case of characteristic
$ 0 $.

\medskip\noindent
{\bf 1.1.~Hypersurface sections.}
Let
$ Y $ be an
$ (m+1) $-dimensional irreducible smooth projective variety
embedded in a projective space
$ \cP $.
For a positive integer
$ i $, let
$ \cO_Y(i) $ be the restriction to
$ Y $ of
$ \cO_{\cP}(i) $ on the projective space
$ \cP $, and define
$$
\cA^{i} = \Gamma(Y,\cO_Y(i)).
$$
Note that the restriction morphism
$$
\Gamma(\cP,\cO_{\cP}(i)) \to \Gamma(Y,\cO_Y(i))
$$
is not surjective in general
(e.g. if
$ Y $ is a hypersurface section of degree
$ i $ of an abelian variety contained in
$ \cP $).
For
$ P \in \cA^{i} $, we denote by
$ X_{P} $ the associated hypersurface section of degree
$ i $ (in a generalized sense unless
$ P $ belongs to the image of the above morphism).
Let
$$
\cV^{i} = \{P\in \cA^{i}\setminus \{0\} : X_{P}\,\,
\text{is smooth}\}.
$$
It is identified with a smooth variety, and
$ \cV^{i}/\bC^{*} $ parametrizes the smooth hypersurface
sections of
$ Y $ of degree
$ i $.

Let
$ Z $ be a closed subvariety of
$ Y $, and put
$$
\cA_{Z}^{i} = \{P \in \cA^{i} : Z \subset X_{P}\},\quad
\cV_{Z}^{i} = \cA_{Z}^{i} \cap \cV^{i}.
$$
In this section we do not necessarily assume condition (0.1)
in Introduction.
Let
$ \delta $ be a positive integer such that
$ \cI_{Z}(\delta) $ is generated by
$ \cA_{Z}^{\delta} $.
It is shown by Kleiman and Altman [13] that condition (0.1)
in Introduction implies the existence of a smooth
hypersurface of degree
$ j $ containing
$ Z $ for any
$ j > \delta' $, where the definition of
$ \delta' $ in loc.~cit. uses the ideal of
$ Z $ in the projective space
$ \cP $ instead of
$ \cI_{Z} $
(in particular
$ \delta \le \delta' $, and we may have a strict inequality
if the degrees of the defining equations of
$ Y $ in
$ \cP $ are bigger than those for
$ Z $ in
$ Y $).
More precisely, they showed the following theorem for
$ d_{i} \ge \delta' + 1 $ and
$ V_{i} = \cA_{Z}^{d_{i}} $,
except possibly for the estimate of
$ \dim\Sing\,X_{\boP} $, which we need later.
Concerning the connectedness, we will need it
for the proof of Theorem (0.4) only in the cases where
$ a \ge 0 $ and either
$ r = 1 $ or
$ 2 $.
Note that the argument is simpler in these cases because
$ X_{P_{1}} $ is smooth.

In Theorem (1.2) below, we will take vector subspaces
$ V_{i} $ of
$ \cA_{Z}^{d_{i}} $, and consider the following
condition:

\begin{itemize}
\item[$(C_{i})$]
$ d_{i} > \delta $ and there exist a vector subspace
$ V'_{i} $ of
$ \cA_{Z}^{\delta} $ generating
$ \cO_{Y}(\delta) $ outside
$ Z $ and a vector subspace
$ V''_{i} $ of
$ \cA^{d_{i}-\delta} $ giving an embedding of
$ Y $ into a projective space such that
$ V_{i} $ is the image of
$ V'_{i} \otimes V''_{i} $.
\end{itemize}

\medskip\noindent
{\bf 1.2.~Theorem} ([13], [19]).
{\it Let
$ a = \dim Y - 1 - \max\{\dim Z_{\{j\}}+j\} $ so that
condition {\rm (0.1)} is satisfied if and only if
$ a \ge 0 $.
Let
$ d_{1} \le \cdots\le d_{r} $ be integers with
$ r \le c_{Z} := \codim_{Y}Z $.
Assume
$ d_{i} \ge \delta $, and let
$ V_{i} $ be a vector subspace of
$ \cA_{Z}^{d_{i}} $ generating
$ \cI_{Z}\otimes\cO_{Y}(d_{i}) $.
Take a general
$ \boP = (P_{1},\dots,P_{r}) $ in
$ V_{1}\mtims\cdots\mtims V_{r} $, and put
$ X_{\boP} = \bigcap_{1\le i\le r}X_{P_{i}} $.
Then
$ X_{\boP} $ is a complete intersection,
$ \dim\Sing\,X_{\boP} \le r-a-2 $, and
$ X_{\boP} \setminus Z $ is smooth.

Assume furthermore that
$ \dim X_{\boP} \ge 1 $, the above condition
$ (C_{i}) $ is satisfied for any integer
$ i $ in
$ [c_{Z},r] $ in the case where
$ \dim Y \le 2\dim Z + a + 2 $ or
$ Z $ is reduced on the complement of a closed subvariety
of dimension
$ < \dim Z $, and
$ (C_{i}) $ is satisfied for any integer
$ i $ in
$ [c_{Z} - 1,r] $ otherwise.
Then $ X_{\boP} \setminus Z $ is connected.
Here
$ [c_{Z},r] = \emptyset $ if
$ r < c_{Z} $.
}

\medskip\noindent
{\it Proof.}
Let
$ \Pi_{r} = V_{1}\mtims\cdots\mtims V_{r} $, and
$$
\Sigma_{r} = \{(\boP,x)\in \Pi_{r}\mtims Z :
\rank\,(dP_{1},\dots,dP_{r})_{x} < r\},
$$
where the
$ dP_{i} $ are defined by taking a local trivialization of
$ \cO(1) $, and the rank is independent of the trivialization
because the
$ X_{P_{i}} $ contain
$ Z $.
We have
$$
T_{x}^{*}Z = m_{Y,x}/(\cI_{Z,x}+m_{Y,x}^{2}),
$$
(where
$ m_{Y,x} $ denotes the maximal ideal), and the
$ P \in V_{i} $ generate
$ (\cI_{Z,x}+m_{Y,x}^{2})/m_{Y,x}^{2} $ for any
$ i $ by the definition of
$ V_{i} $ (taking a local
trivialization of
$ \cO(1) $).
Let
$$
\Sigma_{r,x} = \Sigma_{r}\cap (\Pi_{r}\mtims\{x\}).
$$
By considering the fiber of the projection
$ \Pi_{r}\to\Pi_{r-1} $ and by induction on
$ r $ we get for any
$ x \in Z_{\{\dim Y-j\}} $,
$$
\codim_{\Pi_{r}}\Sigma_{r,x} = \max(j-r+1,0) \ge j-r+1.
$$
By definition of $ a $ we have
$
\codim_{Y}Z_{\{\dim Y-j\}} \ge \dim Y-j+a+1,
$
hence
$$
\codim_{\Pi_{r}\mtims Y}\Sigma_{r} \ge \dim Y - r + a + 2.
$$
For a general
$ \boP \in \Pi_{r} $, this implies
$$
\codim_{Y}\Sigma_{r}\cap (\{\boP\}\mtims Z) \ge \dim Y-r+a+2.
$$
So it remains to show the assertion on the smoothness and
the connectedness.

For the smoothness, we see that a subvariety of
$ (Y\setminus Z)\mtims \Pi_{r} $ defined by the relation
$ x\in X_{\boP} $ for
$ (x,\boP) \in (Y\setminus Z)\mtims \Pi_{r} $ is smooth.
(Indeed, in the case
$ r = 1 $, the variety is defined by
$ \sum_{i} t_{i}P_{i} $ where
$ (P_{i}) $ is a basis of
$ V_{1} $ as a vector space, and
$ (t_{i}) $ is the corresponding coordinate system of
$ V_{1} $.
Furthermore, for any
$ x \in Y\setminus Z $, some
$ P_{i} $ does not vanish on a neighborhood of
$ x $, where the above equation can be divided by this
$ P_{i} $.
The argument is similar for
$ r > 1 $.)
Then the smoothness follows from the Bertini theorem
in characteristic
$ 0 $.

For the connectedness, we proceed by induction on
$ r $.
Let
$ Y' = X_{\boP'} $ for a general
$ \boP'\in \Pi_{r-1} $, where
$ Y' = Y $ if
$ r = 1 $.
Set
$ b = 0 $ if
$ \dim Y \le 2\dim Z + a + 2 $ or
$ Z $ is reduced on the complement of a closed subvariety
of dimension
$ < \dim Z $, and
$ b = 1 $ otherwise.
Assume first
$ r < \codim_{Y}Z - b $.
Since
$ Y'\setminus X_{\boP} $ is a smooth affine variety of
dimension
$ \ge 2 $, its first cohomology with compact supports
vanishes by the weak Lefschetz theorem together with
Poincar\'e duality.
So
$ X_{\boP} $ is connected.
By Hartshorne's connectedness theorem (see (1.3) below)
the connectedness of
$ X_{\boP}\setminus Z $ is then reduced to that
$ \codim_{X_{\boP}}\Sing\,X_{\boP} > 1 $.

Since
$ \Sing\,X_{\boP}\subset Z $, the last condition
is trivially satisfied in the case
$ r < \codim_{Y}Z - 1 $.
If
$ r = \codim_{Y}Z - 1 $ and
$ Z $ is generically reduced as above so that
$ b = 0 $, then
$ X_{\boP} $ is smooth at a general point of
$ Z $ because each
$ V_{i} $ generates
$ \cI_{Z}\otimes \cO_{Y}(d_{i}) $.
Thus the assertion is proved.
Similarly, if
$ r = \codim_{Y}Z - 1 $ and
$ \dim Y \le 2\dim Z + a + 2 $, then
$ \dim\Sing\,X_{\boP} \le r-a-2 < \dim Z $ and
the assertion follows.

Assume now
$ r \ge \codim_{Y}Z - b $.
Consider a rational map of
$ Y' $ to a projective space defined by the restriction of
the linear system
$ V_{r} $ to
$ Y' $.
It induces an embedding of
$ Y' \setminus Z $ because
$ V'_{r}|_{Y'} $ generates
$ \cO_{Y'}(\delta) $ outside
$ Z $ and
$ V''_{r}|_{Y'} $ induces an embedding.
So
$ X_{\boP} \setminus Z $ is isomorphic to a general hyperplane
section of
$ Y' \setminus Z $ for this embedding.
Since
$ Y' \setminus Z $ is a smooth connected variety of dimension
at least
$ 2 $, the Bertini theorem in characteristic
$ 0 $ implies that
$ X_{\boP} \setminus Z $ is connected.
This completes the proof of Theorem (1.2).

\medskip\noindent
{\bf 1.3.~Complement to the proof of Theorem (1.2).}
A local complete intersection is irreducible if it is connected
and the singular locus has codimension
$ > 1 $.
This follows from Hartshorne's connectedness theorem
(see e.g. [9], Th.~18.12),
because a local complete intersection is Cohen-Macaulay.
It also follows from the theory of perverse sheaves [1],
considering the long exact sequence of perverse cohomology
sheaves associated to the distinguished triangle
$ \bQ_{X_{0}}\to\bQ_{X_{1}}\oplus\bQ_{X_{2}}\to
\bQ_{X_{3}}\to $ where
$ X_{0} = X_{1}\cup X_{2} $,
$ X_{3} = X_{1}\cap X_{2} $.
Indeed, if
$ X_{0} $ is a local complete intersection of dimension
$ n $ and
$ X_{3} $ has dimension
$ \le n - 2 $, then
$ \bQ_{X_{0}}[n] $ is a perverse sheaf (i.e.
$ {}^{p}\cH^{n}\bQ_{X_{0}} = \bQ_{X_{0}}[n] $,
see e.g. [7], Th. 5.1.19.) and
$ {}^{p}\cH^{j}\bQ_{X_{3}} = 0 $ for
$ j \ge n-1 $, see [1].
So
$ \bQ_{X_{0}}[n] $ is the direct sum of
$ {}^{p}\cH^{n}\bQ_{X_{i}} $ for
$ i = 1, 2 $, and it is a contradiction if
$ X_{0} $ is connected.

\medskip
The following gives a converse of [13] (and we
reproduce the arguments here for the convenience of the
reader).

\medskip\noindent
{\bf 1.4.~Theorem} [19].
{\it Assume
$ d \ge \delta+1 $.
If there exists a smooth hypersurface section of degree
$ d $ containing
$ Z $,
then condition {\rm (0.1)} in Introduction is satisfied.
}

\medskip\noindent
{\it Proof.}
Let
$ U $ be a non empty smooth open subvariety of
$ Z_{\{\dim Y-j\}} $ such that
$ \dim U = \dim Z_{\{\dim Y-j\}} $ and
$ \cO(1) $ is trivialized over
$ U $.
Let
$$
\cE = \Ker(T^{*}Y|_{U} \to \mcup_{x\in U}T_{x}^{*}Z).
$$
It is a vector bundle of rank
$ j $ over
$ U $, and the associated coherent sheaf
$ \cO_{U}(\cE) $ is generated by
$ dP $ for
$ P\in \cA_{Z}^{\delta} $ (where we fix a trivialization
of
$ \cO(1) $ over
$ U $).
So we may assume that
$ \cE $ is trivialized by
$ P_{1},\dots,P_{j} \in \cA_{Z}^{\delta} $.

Since
$ d \ge \delta+1 $, we see that the
$ dP $ for
$ P\in \cA_{Z}^{d} $ generate the
$ 1 $-jets of
$ \cO_{U}(\cE) $, i.e. they generate
$ \cO_{U}(\cE)_{x}\otimes \cO_{U,x}/m_{U,x}^{2} $ over
$ \bC $ for any
$ x \in U $ (because there are smooth hypersurface
sections of degree
$ d - \delta $ which give local coordinates of
$ U $ at
$ x $).

If the assumption of (1.4) is satisfied, then
$ dP $ for a general
$ P\in \cV_{Z}^{d} $ gives a nowhere vanishing section of
$ \cE $ over
$ U $ because
$ X_{P} $ is smooth.
Since condition (0.1) is equivalent to
$ \dim Z_{\{\dim Y-j\}} < j $ for any
$ j $, the assertion is reduced to the following:

\medskip\noindent
{\bf 1.5.~Lemma.}
{\it Let
$ X $ be a smooth variety of dimension
$ \ge j $, and
$ \cE $ be a trivial vector bundle of rank
$ j $ over
$ X $.
Let
$ V $ be a finite dimensional vector subspace of
$ \Gamma(X,\cE) $ which generates the
$ 1 $-jets at every point of
$ X $.
Then for a general
$ \sigma\in V $, the zero locus
$ \sigma^{-1}(0) \subset X $ is non empty.
}

\medskip\noindent
{\it Proof.}
Since
$ \cE $ is trivial, we may identify
$ \sigma\in \Gamma(X,\cE) $ with a morphism of
$ X $ to a vector space
$ E $.
We have a natural morphism
$ \theta : V\mtims X\to E $ which sends
$ (\sigma,x) $ to
$ \sigma(x) $.
Let
$ Y = \theta^{-1}(0) $.
Then it is enough to show that
$ Y $ is dominant over
$ V $.

Since
$ V $ generates the
$ 1 $-jets and
$ \dim X \ge j $, there exists
$ y = (\sigma,x) \in Y $ such that the differential
$ d_{x}\sigma : T_{x}X \to T_{\sigma(x)}E = E $ of
$ \sigma : X \to E $ (by the above identification)
is surjective.
Consider the commutative diagram
$$
\CD
0 @>>> T_{x}X @>>> T_{y}(V\mtims X) @>>> T_{\sigma}V
@>>> 0
\\
@. @VV{d_{x}\sigma}V @VV{d_{y}\theta}V @VVV
\\
0 @>>> E @>>> E @>>> 0 @>>> 0
\endCD
$$
where the top row is induced by the inclusion
$ \{\sigma\}\mtims X\to V\mtims X $ and the projection
$ V\mtims X \to V $.
Then the surjectivity of
$ d_{x}\sigma $ implies that of
$ \Ker\,d_{y}\theta\to T_{\sigma}V $ by the snake lemma.
So the assertion follows because
$ \Ker\,d_{y}\theta = T_{y}Y $.
This completes the proof of Theorem (1.4).

\medskip\noindent
{\bf 1.6.~Voisin's proof of Conjecture (0.2) in the smooth case.}
Consider a rational morphism of
$ Y $ to a projective space
$ \cP $ defined by the linear system
$ \cA_{Z}^{\delta+1} $.
It induces an embedding of
$ Y \setminus Z $, in particular, it is birational to the image.
Let
$ \tY $ be the closure of the image.
If
$ Z $ is smooth, we see that
$ \tY $ is the blow-up of
$ Y $ along
$ Z $.
For a smooth hyperplane section
$ \tY_{s} $ of
$ \cP $ and the corresponding smooth hypersurface section
$ Y_{s} $ of
$ Y $, there is a morphism
$ H^{m}(\tY_{s})^{\van} \to H^{m}(Y_{s})^{\van} $ by
the functoriality of the Gysin morphisms, and it is injective
by the irreducibility of
$ \{H^{m}(\tY_{s})^{\van}\} $.
Calculating the cohomology of the blow-up,
we see that the dimension of its cokernel coincides with
the number of the
$ (m/2) $-dimensional irreducible components of
$ Z $.
So Conjecture (0.2) is proved in this case.

\bigskip\bigskip
\centerline{\bf 2. Topological Method}

\bigskip\noindent
In this section, we prove the nonvanishing of certain
extension classes when the condition
$ d \ge \delta + 2 $ in Theorem (0.4) is satisfied.

\medskip\noindent
{\bf 2.1.~Exact sequences.}
Let
$ Y $ be a connected smooth complex algebraic variety, and
$ X $ be a divisor on
$ Y $ with the inclusion
$ i : X \to Y $.
Put
$ U = Y \setminus X $ with the inclusion
$ j : U \to Y $.
Let
$ f : Y \to S $ be a proper morphism to a smooth variety
$ S $,
and let
$ g = f\scirc i : X \to S, h = f\scirc j : U \to S $ be the
restrictions of
$ f $.
Put
$ m = \dim Y - \dim S $.
We have a long exact sequence of constructible sheaves
$$
\aligned
&\to R^{m-1}f_{*}\bQ_{Y} \buildrel{\gamma_{m-1}}\over\to
R^{m-1}g_{*}\bQ_{X}
\\
&\qquad\qquad\qquad\to R^{m}h_{!}\bQ_{U}\to R^{m}f_{*}\bQ_{Y}
\buildrel{\gamma_{m}}\over\to R^{m}g_{*}\bQ_{X}\to,
\endaligned
\leqno(2.1.1)
$$
using the distinguished triangle
$ j_{!}\bQ_{U} \to \bQ_{Y} \to i_{*}\bQ_{X}
\buildrel{+1}\over\to $ together with the functor
$ \boR f_{*} = \boR f_{!} $.
Since
$ f, g $ are proper, the base change holds so that the stalk of
the direct image is isomorphic to the (relative) cohomology
of the fiber.
Let
$$
\cF = R^{m}h_{!}\bQ_{U},\quad \cF' = \Coker \,\gamma_{m-1},
\quad \cF'' = \Ker \,\gamma_{m},
$$
so that we have a short exact sequence of constructible
sheaves
$$
0 \to \cF' \to \cF \to \cF'' \to 0.
\leqno(2.1.2)
$$
Let
$ X_{s} =g^{-1}(s) $, etc.
If
$ f, g $ are smooth projective and
$ X_{s} $ is a hypersurface section of
$ Y_{s} $ for
$ s \in S $, then we have
$$
\cF'_{s} = H^{m-1}(X_{s},\bQ)^{\van},\quad
\cF_{s} = H^{m}(Y_{s},X_s,\bQ),\quad
\cF''_{s} = H^{m}(Y_{s},\bQ)^{\prim}.
\leqno(2.1.3)
$$
We will assume that they are nonzero (because otherwise
the extension class vanishes).

If
$ f, g $ are smooth morphisms, then (2.1.1) and (2.1.2) are
exact sequences of local systems, which underlie naturally
variation of mixed Hodge structures, see [4], [25]
(and also [22]).
In the application we will also consider the dual of (2.1.2)
$$
0 \to \cF''{}^{*}\to \cF^{*} \to \cF'{}^{*} \to 0,
$$
where
$ {}^{*} $ denotes the dual variation of mixed Hodge structure.
Note that
$$
\cF'{}^{*} = \cF'(m-1),\quad
\cF''{}^{*} = \cF''(m),\quad
\cF^{*}_{s} = H^{m}(Y_{s}\setminus X_{s},\bQ)(m),
$$
where
$ (m) $ denotes the Tate twist, see [4].

\medskip\noindent
{\bf 2.2.~Lefschetz pencils.}
With the above notation, assume
$ S = \bP^{1} $ and
$ g : X \to S $ is a Lefschetz pencil of a smooth
projective variety
$ Y_{0} $.
Let
$ H^{m-1}(X_{s},\bQ)^{\van} $ be the subgroup of
$ H^{m-1}(X_{s},\bQ) $ generated by the vanishing cycles,
and assume it nonzero.
Let
$ s_{0} $ be any point in the discriminant of the Lefschetz
pencil.
Since the discriminant in the parameter space of hypersurfaces
is irreducible,
the last assumption is equivalent to the surjectivity of
the restriction morphism
$$
H^{m-1}(X_{s},\bQ) \to H^{m-1}(B_{s},\bQ)
\leqno(2.2.1)
$$
for
$ s $ sufficiently near
$ s_{0} $, where
$ B $ is a small ball in
$ X $ around the unique singular point of
$ X_{s_{0}} $, and
$ B_{s} := B\cap X_{s} $ is called a Milnor fiber, see [17].
This implies that the cospecialization morphism
$$
H^{m}(X_{s_{0}},\bQ) \to H^{m}(X_{s},\bQ)
\leqno(2.2.2)
$$
is an isomorphism, using a long exact sequence.
So we get
$$
R^{j}g_{*}\bQ_{X}\,\,\text{is a constant sheaf on}\,\,
S \,\,\text{for any}\,\, j \ne m-1,
\leqno(2.2.3)
$$
see [12], where the case
$ j = m $ follows from the above argument,
and the other cases are easy.
As a corollary, we get in this case
$$
\cF'' \,\,\text{in (2.1.2) is a constant sheaf.}
\leqno(2.2.4)
$$

The above argument also implies for
$ j = m - 1 $ that
$ R^{m-1}g_{*}\bQ_{X} $ is a (shifted) intersection complex,
i.e.
$$
R^{m-1}g_{*}\bQ_{X} = j_{*}j^{*}R^{m-1}g_{*}\bQ_{X},
\leqno(2.2.5)
$$
where
$ j : S' \to S $ is the inclusion of a dense open
subvariety over which
$ g $ is smooth.
(This also follows from the local invariant cycle theorem
[3] or the decomposition theorem [1].)

\medskip
The following proposition was proved in [19] using a
generalization of the Picard-Lefschetz formula together with
an assertion concerning the vanishing cycles of a Lefschetz
pencil and related to the classical work of Lefschetz and
Poincar\'e (see [14], [15], [28]).
We give here a simple proof of the proposition using the above
cohomological property of the Lefschetz pencil.

\medskip\noindent
{\bf 2.3.~Proposition} ({Lefschetz pencil case}) [19].
{\it With the notation and the assumptions of {\rm (2.1)},
assume
$ S = \bP^{1} $,
$ Y = Y_{0} \mtims S $,
$ f = \pr_{2} $, and
$ g : X \to S $ is a Lefschetz
pencil of
$ Y_{0} $.
Let
$ S' $ be any non empty open subvariety of
$ S $ over which
$ g $ is smooth.
Then for any nonzero local subsystem
$ \cG $ of
$ \cF''|_{S'} $, the composition of the inclusion
$ \cG \to \cF''|_{S'} $ with the extension class defined by
the restriction of the short exact sequence {\rm (2.1.2)}
to
$ S' $ is nontrivial as an extension of local systems.
}

\medskip\noindent
{\it Proof.}
Since the local system
$ \cF''|_{S'} $ is constant, we may assume that
$ \cG $ has rank
$ 1 $, and is generated by
$ u \in H^{0}(S',\cF''|_{S'}) = H^{m}(Y_{0},\bQ)^{\prim} $.
Assume
$ u $ is the image of
$ v\in H^{0}(S',\cF|_{S'}) $.
Then it gives a section of (2.1.2) on
$ \cG $.
So it is enough to show that
$ u = 0 $ in this case.

We see that
$ \cF $ is a (shifted) intersection complex by (2.1.2),
because
$ \cF' $ and
$ \cF''$ are (shifted) intersection complexes with support
$ S $.
So
$ H^{0}(S,\cF) = H^{0}(S',\cF|_{S'}) $, and we may replace
$ S' $ with
$ S $ or any nonempty open subvariety of
$ S $.
Thus we may assume that
$ S' = S \setminus \{s_{0}\} $ and
$ X_{s_{0}} $ is smooth.

Let
$ h' : U' \to S' $ be the restriction of
$ h $ over
$ S' $, where
$ U' = Y'\setminus X' $.
Consider the Leray spectral sequence
$$
E_{2}^{p,q} = H^{p}(S',R^{q}h'_{!}\bQ_{U'})
\Rightarrow H^{p+q}(Y',X';\bQ).
\leqno(2.3.1)
$$
This degenerates at
$ E_{2} $, because
$ E_{2}^{p,q} = 0 $ unless
$ p = 0 $ or
$ 1 $.
Thus we get
$ w\in H^{m}(Y',X';\bQ) $ whose image in
$ H^{0}(S',\cF|_{S'}) $ is
$ v $.
Its image in
$ H^{0}(S',\cF''|_{S'}) = H^{m}(Y',\bQ) = H^{m}(Y_0,\bQ) $ is
$ u $, where the last isomorphism follows from
$ Y' = Y_{0}\mtims \bA^{1} $.
Then the image of
$ u $ in
$ H^{m}(X',\bQ) $ vanishes.
But this is induced by the restriction morphism under
the birational morphism
$ X'\to Y_{0} $.
So we can verify that
$ u $ belongs to the image of the Gysin morphism under the
inclusion
$ X_{s_{0}}\to Y_{0} $, and we get
$ u = 0 $ because
$ u $ is primitive.
This completes the proof of Proposition (2.3).

\medskip
We also give an outline of the original proof of Proposition
(2.3).
We start with the explanation of a generalized Picard-Lefschetz
formula.

\medskip\noindent
{\bf 2.4.~Generalized Picard-Lefschetz formula.}
Let
$ \cF $ be a constructible sheaf on a curve
$ S $ with a local coordinate
$ t $.
Let
$ \psi_{t}\cF, \varphi_{t}\cF $ denote the nearby and
vanishing cycles, see [5].
Then we have natural morphisms
$$
\can : \psi_{t}\cF \to \varphi_{t}\cF,\quad
\var : \varphi_{t}\cF \to \psi_{t}\cF,
$$
such that
$$
T - id = \var\scirc \can : \psi_{t}\cF \to \psi_{t}\cF,
\leqno(2.4.1)
$$
where
$ T $ is the monodromy.
It is well known that the functors
$ \psi_{t}, \varphi_{t} $ commute with the functor
assigning the dual, and duality exchanges
$ \can $ and
$ \var $ up to a sign, see e.g. [11], [21].
In the case of a Lefschetz pencil, we can identify the morphism
$ \can $ with the restriction to the Milnor fiber (2.2.1).

Assume that
$ \varphi_{t}\cF \simeq \bQ $, and
$ \can $,
$ \var $ are nonzero.
Let
$ \gamma_{0} $ be a generator of
$ \varphi_{t}\cF $, and
$ \gamma $ be its image in
$ \psi_{t}\cF $ by
$ \var $.
Let
$ \cF^{*} $ be the (shifted) dual of
$ \cF $ which is defined by
$ \boR\cHom(\cF,\bQ_{S}) $.
Let
$ \gamma_{0}^{*} $ be the generator of
$ \varphi_{t}\cF^{*} $ such that
$ \langle \gamma_{0}^{*},\gamma_{0} \rangle = 1 $.
Let
$ \gamma^{*} $ be its image in
$ \psi_{t}\cF^{*} $ by
$ \var $.
Then we have a generalized Picard-Lefschetz formula
$$
T(u) - u = \pm \langle \gamma^{*},u\rangle\gamma\quad
\text{for}\,\,u \in \psi_{t}\cF,
\leqno(2.4.1)
$$
because
$ \langle \gamma_{0}^{*}, \can(u)\rangle =
\pm \langle \var(\gamma_{0}^{*}),u\rangle $, see also [6].
This was proved in [19] for the cohomology of the complement
of a hypersurface section.

\medskip\noindent
{\bf 2.5.~Outline of the original proof of Proposition (2.3)}
(see [19]).
It is enough to show that
$ \cF|_{S'} $ has no global sections.
In this case, the stalk of
$ \cF^{*}|_{S'} $ is
$ H_{m}(Y_{s},X_{s}) $, and
$ \gamma^{*} $ can be constructed explicitly using the
ball and the Milnor fiber around the critical point (and
this coincides with the construction in [14], [15], [28]),
because we can identify the morphism
$ \can $ with the restriction to the
complement of the Milnor fiber in the ball, see [19].
Furthermore, considering
$ \gamma^{*} $ at any points of the discriminant of the
Lefschetz pencil, they generate
$ H_{m}(Y_{s},X_{s}) $.
(This is closely related to the classical work of
Lefschetz and Poincar\'e, and seems to have been known to
some people, see [14] and also [15], [28].)
So the local system
$ \cF|_{S'} $ has no global section,
and the assertion follows.

\medskip\noindent
{\bf 2.6.~Non Lefschetz fibration case.}
Let
$ Y $ be a connected smooth projective variety embedded in
a projective space
$ \cP $, and
$ X $ be a hypersurface section of
$ Y $ with at most isolated singularities.
We assume
$ m := \dim X = \dim Y - 1 \ge 1 $.
Let
$ Z $ be an irreducible component of
$ X $
(hence $ Z = X $ if $ m > 1 $).
Let
$ d $ be an integer
$ \ge 2 $, and
$ S $ be the parameter space of hypersurfaces of degree
$ d $ of
$ \cP $ whose intersections
$ Z_{s}, Y_{s} $ with
$ Z, Y $ are smooth divisors on
$ Z \setminus \Sing\, X $ and
$ Y $ respectively
(in particular,
the hypersurfaces parametrized by
$ S $ do not meet
$ \Sing\, X $).

Let
$ \{H^{m-1}(Z_{s})_{Z}^{\van}\}_{s\in S} $ be the local
subsystem of
$ \{H^{m-1}(Z_{s})\}_{s\in S} $ generated by the
vanishing cycles at general points of the
discriminant of the morphism
$ \bigsqcup_{s\in \oS}Z_{s}\to \oS $, where
$ \oS $ is the parameter space of all the hypersurfaces
of degree
$ d $ of
$ \cP $, and
$ \bigsqcup_{s\in \oS}Z_{s} $ denotes the total space of
the associated family of hypersurfaces.
If
$ m > 1 $ and
$ X = Z $, let
$ \{H^{m-1}(Z_{s})^{\van}\}_{s\in S} $ be the local subsystem
generated by the vanishing cycles for the inclusion
$ Z_{s} \to Y_{s} $.
By the Picard-Lefschetz formula, the latter is the orthogonal
complement of the injective image of the cohomology of
$ Y_{s} $ (or
$ Y $ using the weak Lefschetz theorem), and hence contains
the former.
If
$ X = Z $ is smooth, they coincide because they are the
orthogonal complement of the injective image of
$ H^{m-1}(Y) $.
If
$ m = 1 $, let
$ H^{m-1}(Z_{s})^{\van} = H^{m-1}(Z_{s})_{Z}^{\van}\,
(= \tH^{m-1}(Z_{s})) $.

Let
$ S' $ be a dense open subvariety of
$ S $, and
$ \tL $ be any local system on
$ S' $ such that
$$
\{H^{m-1}(Z_{s})_{Z}^{\van}\}_{s\in S'} \subset \tL \subset
\{H^{m-1}(Z_{s})^{\van}\}_{s\in S'},
$$
and the restriction of the intersection pairing to
$ \tL $ is nondegenerate.
Let
$ \tL^{\bot} $ be the orthogonal complement of
$ \tL $ in
$ \{H^{m-1}(Z_{s})\} $.
Note that the restrictions of the intersection pairing to the
injective image of
$ H^{m-1}(Y) $ and to
$ H^{m-1}(Z_{s})^{\van} $ are non degenerate using
Hodge theory (or [1] because it is essentially equivalent to
the hard Lefschetz theorem).

Consider the kernel of the composition
$$
H^{m}(Y_{s},Z_{s}) \to H^{m}(Y_{s})\to H^{m+2}(Y)(1),
$$
where the last morphism is the Gysin morphism.
Let
$ H^{m}(Y_{s},Z_{s})_{\tL}^{\van} $ be the quotient of the
kernel by the image of
$ \tL^{\bot} $.

\medskip
The following is a generalization of a construction in
the surface case in [19], and gives a topological proof of
variants of Propositions (3.2) in the non Lefschetz case
and (3.3) in the surface case.

\medskip\noindent
{\bf 2.7.~Proposition} ({Non Lefschetz fibration case}).
{\it With the above notation and assumptions,
we have a short exact sequence of local
systems on
$ S' $
$$
0 \to \tL \to \{H^{m}(Y_{s},Z_{s})_{\tL}^{\van}\} \to
\{H^{m}(Y_{s})^{\van}\}\to 0,
\leqno(2.7.1)
$$
and it does not split if the first and last terms are nonzero.
}

\medskip\noindent
{\it Proof.}
The exactness of (2.7.1) is clear by definition.
To show the non splitting of (2.7.1), we may assume
$ S' = S $ using the direct image by
$ S ' \to S $.
We take a smooth point
$ O $ of
$ X $ contained in
$ Z $, and consider a hypersurface
$ H_{0} $ in the ambient projective space which intersects
$ X, Y $ transversely at smooth points outside
$ O $ and whose intersections
with
$ Z $,
$ Y $ have an ordinary double point at
$ O $.
(Such a hypersurface exists because the degree of the
hypersurface is at least
$ 2 $.)
Then we consider a generic deformation
$ \{H_{\lambda}\} $ of
$ H_{0} $, parametrized by
$ \Lambda := \Delta^{m+2} $ where
$ \Delta $ is a sufficiently small open disk.
Using the embedding by
$ \cO_{Y}(d) $ (which changes hypersurfaces into hyperplanes)
and an appropriate projection between projective spaces,
we may assume locally
$$
\aligned
Y
&= \{x_{m+2} = \msum_{1\le i \le m+1}x_{i}^{2} +
\Psi(x_{1},\dots,x_{m+1})\},
\\
Z
&= \{x_{m+1} = 0\} \subset Y,
\\
H_{\lambda}
&= \{x_{m+2} = \msum_{1\le i\le m+1}2a_{i}x_{i} + a_{m+2}\}\quad
\text{for}\,\,
\lambda = (a_{i}) \in \Delta^{m+2},
\endaligned
$$
in a
$ (m+2) $-dimensional projective space, where
$ (x_{i}) $ is an affine coordinate system and
$ \Psi(x_{1},\dots,x_{m+1})\in (x_{1},\dots,x_{m+1})^{3} $.

Let
$ Z_{\lambda} = Z\cap H_{\lambda} $,
$ Y_{\lambda} = Y\cap H_{\lambda} $.
The discriminant of
$ \{Z_{\lambda}\}_{\lambda\in \Lambda} $ and
$ \{Y_{\lambda}\}_{\lambda\in \Lambda} $ in
$ \Lambda $ are defined
respectively by
$$
\aligned
D_{Z}
&= \{a_{m+2} + \msum_{1\le i \le m}a_{i}^{2} +
\Phi_{Z}(a_{1},\dots,a_{m}) = 0\},
\\
D_{Y}
&= \{a_{m+2} + \msum_{1\le i\le m+1}a_{i}^{2} +
\Phi_{Y}(a_{1},\dots,a_{m+1}) = 0\},
\endaligned
$$
where
$ \Phi_{Z}(a_{1},\dots,a_{m}) \in
(a_{1},\dots,a_{m})^{3} $,
$ \Phi_{Y}(a_{1},\dots,a_{m+1}) \in
(a_{1},\dots,a_{m+1})^{3} $.

Take
$ c \in \Delta\setminus \{0\} $, and put
$ 0' = (0)\in\Delta^{m} $.
Define
$$
\aligned
\Lambda_{c}
&= \{0'\}\mtims\{c\}\mtims\Delta,\quad
\{\alpha\} = \Lambda_{c}\cap D_{Z},\quad
\{\beta\} = \Lambda_{c}\cap D_{Y},
\\
\Lambda'_{c}
&=\Lambda_{c}\setminus\{\alpha,\beta\}.
\endaligned
$$
Let
$ B $ be a sufficiently small open ball around
$ O $ in the ambient space, and assume
$ \Delta $ is much smaller.
Let
$$
Z_{B} = Z\cap B,\quad Z_{B,\lambda} = Z_{\lambda}\cap B
\,\,\, \text{(similarly for}\,\, Y).
$$
For
$ \lambda\in\Lambda_{c}\setminus\{\beta\} $,
it is well known that
$$
\aligned
&\tH^{j}(Y_{B,\lambda}) = \bQ\quad
\text{for}\,\,j = m,\,\,
\text{and}\,\,0\,\,\text{otherwise},
\\
&\tH^{j}(Y_{B,\beta}) = 0\quad
\text{for any}\,\, j.
\endaligned
\leqno(2.7.2)
$$
and similarly for
$ Z_{B,\lambda} $ with
$ \beta $ replaced by
$ \alpha $, and
$ m $ by
$ m - 1 $.

The local monodromy of
$ H^{m-1}(Z_{B,\lambda}) $ (resp.
$ H^{m}(Y_{B,\lambda}) $) around
$ \alpha $ (resp.
$ \beta $) is
$ (-1)^{m} $ (resp.
$ (-1)^{m+1} $).
By (2.7.2) we have the short exact sequences

\bigskip\noindent
(2.7.3)\qquad\qquad
$ 0\to H^{m-1}(Z_{B,\lambda})\to
H^{m}(Y_{B,\lambda},Z_{B,\lambda})\to
H^{m}(Y_{B,\lambda})\to 0 $,

\bigskip\noindent
(2.7.4)\qquad\qquad
$ 0\to H_{c}^{m-1}(Z_{B,\lambda})\to
H_{c}^{m}(Y_{B,\lambda}\setminus Z_{B,\lambda})\to
H_{c}^{m}(Y_{B,\lambda})\to 0 $.

\bigskip\noindent
By definition, we have an injective morphism of (2.7.4)
to (2.7.1).
Taking the dual, we get a surjective morphism of (2.7.1) to
(2.7.3) because the intersection form is non degenerate on
$ \tL $.
Let
$ e, e' $ denote the extension classes defined by
(2.7.3) and (2.7.4) respectively.
Then the assertion is reduced to Lemma (2.8) below.
Indeed, if
$ m $ is odd and the composition of
$ e $ with the surjection
$$
\{H^{m}(Y_{\lambda})^{\van}\}_{\lambda\in\Lambda'_{c}} \to
\{H^{m}(Y_{B,\lambda})\}_{\lambda\in\Lambda'_{c}}
$$
vanishes, then
$ \{H^{m}(Y_{B,\lambda},Z_{B,\lambda})\} $ is isomorphic to
a subquotient of
$ \{H^{m}(Y_{\lambda})^{\van}\} $ (by the nontriviality of
$ e $) and there is a surjection
from the kernel of the above morphism to
$ \{H^{m-1}(Z_{B,\lambda})\} $, but this contradicts
the nontriviality of the monodromy of
$ \{H^{m-1}(Z_{B,\lambda})\} $ around
$ \alpha $.
The argument is similar for
$ m $ even.

\medskip\noindent
{\bf 2.8.~Lemma.}
{\it The extension classes
$ e, e' $ defined by {\rm (2.7.3), (2.7.4)}
do not vanish if
$ m $ is odd or even respectively.
}

\medskip\noindent
{\it Proof.}
We first show the non splitting of (2.7.3) for
$ m $ odd.
In this case, the monodromy of
$ \{H^{m}(Y_{B,\lambda},Z_{B,\lambda})\} $ around
$ \beta $ is unipotent, and we have to show that it is not
the identity.
Let
$ D\subset\Delta_{c} $ be a sufficiently small open disk
with center
$ \beta $, and restrict the local systems to
$ D^{*} := D\setminus \{\beta\} $.
The short exact sequence of local systems (2.7.3) on
$ D^{*} $ is naturally extended to
$ D $ as constructible sheaves,
by using the higher direct images as in (2.1)
instead of cohomology.

More precisely, let
$ Z_{B,D}, Y_{B,D} $ be the restriction of
$ Z_{B}, Y_{B} $ over
$ D $.
Let
$ U = Y_{B,D}\setminus Z_{B,D} $ with the inclusion
$ j_{U} : U \to Y_{B,D} $
Then the short exact sequence of local systems on
$ D^{*} $ defined by (2.7.3) is
naturally extended to a short exact sequence of
constructible sheaves on
$ D $
$$
0 \to \cG' \to \cG \to \cG'' \to 0.
\leqno(2.8.1)
$$
This is also defined by the higher direct images of
$ (j_{U})_{!}\bQ_{U} $, etc. under
$ Y_{B,D} \to D $.
Here the base change holds, because
$ (Y_{\lambda}, Z_{\lambda}) $ is transversal to the boundary
of
$ B $.
In particular, the stalk at
$ \beta $ is given by (2.7.3).

Using (2.7.2), we see that
$ \cG' = \bQ_{D} $,
$ \cG'' = (j_{D})_{!}\bQ_{D^{*}} $, where
$ j_{D} : D^{*} \to D $ is the inclusion morphism.
Let
$ t $ be a local coordinate around
$ \beta $.
Then
$ \varphi_{t}\cG' = 0 $ and
$ \psi_{t}\cG' = \psi_{t}\cG'' = \varphi_{t}\cG'' = \bQ $ using
the bijectivity of
$ \can $ (see (2.4)) for
$ \cG'' $.
The last assertion follows from the distinguished triangle
of functors
$$
i_{0}^{*}\to\psi_{t}\buildrel{\can}\over\longrightarrow
\varphi_{t} \to,
\leqno(2.8.2)
$$
where
$ i_{0} : \{\beta\} \to D $ denotes an inclusion, see [5].
Since
$ \psi_{t}, \varphi_{t} $ are exact functors, we get
$ \rank\,\psi_{t}\cG = 2 $,
$ \rank\,\varphi_{t}\cG = 1 $ together with the surjectivity of
$ \can : \psi_{t}\cG \to \varphi_{t}\cG $.

By (2.4.1) the assertion is thus reduced to the injectivity of
$ \var : \varphi_{t}\cG \to \psi_{t}\cG $.
By duality (see (2.4)) it is equivalent to the surjectivity of
$ \can : \psi_{t}\cG^{*} \to \varphi_{t}\cG^{*} $,
where
$ \cG^{*} $ is the (shifted) dual of
$ \cG $ which is defined by
$ \boR\cHom(\cG, \bQ_{D}) $.

By duality [27]
$ \cG^{*} $ is isomorphic to the derived direct image with
compact supports of
$ \boR (j_{U})_{*}\bQ_{U} $.
The stalk at
$ \beta $ of the higher direct image sheaves is
$$
H_{c}^{j}(Y_{B,\beta},\boR (j_{\beta})_{*}
\bQ_{Y_{\beta}\setminus Z_{\beta}}),
\leqno(2.8.3)
$$
(where
$ j_{\beta} : Y_{\beta}\setminus Z_{\beta}\to Y_{\beta} $
denotes a natural inclusion), because
$ Z_{\beta} $ and
$ Y_{\beta} $ are smooth on a neighborhood of
$ Z_{\beta} $.
So the assertion is reduced to the vanishing of (2.8.3) for
$ j = m + 1 $ by the distinguished triangle (2.8.2)
because it gives the cokernel of
$ \can $.

Using a one-parameter deformation
$ t\Psi \,(t\in \bC) $, we can reduce the assertion to the case
$ \Psi = 0 $.
So we may assume that
$ Y_{B,\beta} $ is the intersection of
$ B $ with an affine cone
$ Y_{\beta} $ of a nonsingular conic, and
$ Z_{B,\beta} $ is its intersection with a linear space passing
near the origin.
Then we may omit the subscript
$ B $ by replacing
$ B $ with arbitrary large open balls and taking the limit.
So the assertion follows from Artin's generalization of the weak
Lefschetz theorem [1].
Indeed, the closure
$ \oZ_{\beta} $ of
$ Z_{\beta} $ in the projective space intersects transversely
the divisor at infinity
$ \oY_{\beta}\setminus Y_{\beta} $ of
$ Y_{\beta} $ so that we can first take the direct image with
compact supports by
$ Y_{\beta}\setminus Z_{\beta}\to
\oY_{\beta}\setminus \oZ_{\beta} $, and then the usual direct
image by
$ \oj_{\beta} : \oY_{\beta}\setminus \oZ_{\beta}\to
\oY_{\beta} $.
Thus the assertion is proved for
$ m $ odd.

In the case
$ m $ is even, let
$ D $ be a sufficiently small neighborhood of
$ \alpha $ in
$ \Delta_{c} $.
For a variety
$ V $ of pure dimension
$ r $, let
$ \tH_{c}^{j}(V) = H_{c}^{j}(V) $ for
$ j \ne 2r $, and define
$ \tH_{c}^{2r}(V) $ to be the kernel of
$ \Tr : H_{c}^{2r}(V) \to \bQ(-r) $ which is induced by the
canonical morphism of
$ \bQ_{V}(r)[2r] $ to the dualizing complex
$ \bD_{V} $.
Then
$ \tH_{c}^{j}(Z_{B,\alpha}) = \bQ $ for
$ j = m-1, m $ by calculating the cohomology of
$ \partial Z_{B,\alpha} $ and using the long exact sequence
$ \to H^{j-1}(\partial Z_{B,\alpha})
\to H_{c}^{j}(Z_{B,\alpha})
\to H^{j}(Z_{B,\alpha}) \to $.
In particular, (2.7.4) for
$ \lambda = \alpha $ is not exact, and is extended to
a long exact sequence.
So we have to consider a distinguished triangle in
the derived category of sheaves on
$ D $,
which is defined by using the direct images with compact
supports under the morphisms of
$ Z_{B,D}, Y_{B,D}\setminus Z_{B,D}, Y_{B,D} $ to
$ D $ (in particular, the base change holds).
They contain some shifted constant sheaves
which are annihilated by taking the
reduced cohomology in (2.7.2).
So deleting them, we get a short exact sequence of
shifted perverse sheaves as in (2.8.1)
such that the stalk of the
$ 0 $-th cohomology sheaf is given by (2.7.4).
By the above calculation, we have
$ \cG' = \boR (j_{D})_{*}\bQ_{D^{*}} $ and
$ \cG'' = \bQ_{D} $ in this case.

Using an argument similar to the case
$ m $ odd, we see that
$ \var : \varphi_{t}\cG \to \psi_{t}\cG $ is injective
(reducing to the assertion for
$ \cG' $ in this case)
where
$ t $ is a coordinate around
$ \alpha $.
So it is enough to show the nontriviality of
$ \can : \psi_{t}\cG \to \varphi_{t}\cG $.
By the triangle (2.8.2) this is equivalent to
$ H_{c}^{m}(Y_{B,\alpha}\setminus Z_{B,\alpha}) = \bQ $,
because it gives the kernel of
$ \can $.
By duality it is further equivalent to
$$
H^{m}(Y_{B,\alpha}\setminus Z_{B,\alpha}) = \bQ.
\leqno(2.8.4)
$$

By the same argument as above, we may omit the subscript
$ B $ by assuming that
$ Y_{\alpha} $ is the restriction of a nonsingular conic
$ \oY_{\alpha} $ in
$ \bP^{m+1} $ to the affine space
$ \bA^{m+1} $,
$ Z_{\alpha} $ is its intersection with a hyperplane
$ H $, which is an affine cone of a nonsingular conic in
$ \bP^{m-1} = H\cap\bP^{m} $, and the divisor at infinity
$ \partial Y_{\alpha} = \oY_{\alpha}\setminus Y_{\alpha} $
of
$ Y_{\alpha} $ is smooth.
Then, using a projection from the vertex of the affine cone
$ Z_{\alpha} $, we see that
$ Y_{\alpha}\setminus Z_{\alpha} $ is isomorphic to
the complement of the union of
$ \partial Y_{\alpha} $ and
$ H\cap \bP^{m} $ in
$ \bP^{m} = \bP^{m+1}\setminus \bA^{m+1} $.
So (2.8.4) follows considering the long
exact sequence containing the Gysin morphism by
$ \partial Y_{\alpha}\setminus H \to \bA^{m} =
\bP^{m}\setminus H $.
This completes the proofs of Lemma (2.8) and Proposition
(2.7).

\bigskip\bigskip
\centerline{\bf 3. Hodge-Theoretic Method}

\bigskip\noindent
In this section, we prove the nonvanishing of certain
extension classes when the assumption on a differential
form in Theorem (0.4) is satisfied.

\medskip\noindent
{\bf 3.1.~Extension groups.}
With the notation and the assumptions of (2.1), we
assume in this section that
$ f, g $ are smooth (by restricting
$ S $), and
$ \cF', \cF'' $ are nonzero.
We will consider whether (2.1.2) splits in the category of
local systems.
Let
$$
\cH = \cHom(\cF'',\cF').
$$
Then there is a canonical isomorphism
$$
\Ext^{1}(\bQ_{S},\cH) = \Ext^{1}(\cF'',\cF'),
\leqno(3.1.1)
$$
where the extension group is taken in the category
of admissible variations of mixed Hodge structures [25]
(or equivalently, in that of mixed Hodge modules [22]).
These groups are identified (using [2]) with the scalar
extension of the group of admissible normal functions,
which are sections of the family of Jacobians
$ \{J(\cH_{s})\}_{s\in S} $ satisfying some good conditions
[23].
(Here
$ \cH $ also denotes a variation of mixed
$ \bZ $-variation of Hodge structure whose scalar extension
is
$ \cH $).
We have furthermore a short exact sequence
$$
0 \to \Ext^{1}(\bQ,H^{0}(S,\cH)) \to
\Ext^{1}(\bQ_{S},\cH) \buildrel{r}\over\to
\Hom(\bQ,H^{1}(S,\cH)) \to 0,
\leqno(3.1.2)
$$
see [22], [29], where
$ \Ext $ and
$ \Hom $ are taken in the category of mixed Hodge structures
or that of admissible variations of mixed Hodge structures.
(Using [22], this follows from the adjoint relation between
the direct image and the pull-back of mixed Hodge modules by
$ S \to pt $.
Using [29], we get the above short exact sequence with
the cohomology in the last term replaced by the intersection
cohomology, which is a subgroup of the cohomology in this
case.
But this is enough for our purpose, although we can show that
these two give the same by taking
$ \Hom $, see also [8].)

Let
$ \cF_{\bQ},_{ }\cH_{\bQ}, $ etc.
denote the underlying local systems.
Then we have
$$
\Ext^{1}(\cF''_{\bQ},\cF'_{\bQ}) = H^{1}(S,\cH_{\bQ}).
\leqno(3.1.3)
$$
This is compatible with the last morphism
$ r $ in (3.1.2).
Let
$ e $ be the extension class defined by the short exact
sequence (2.1.2).
Then (2.1.2) splits in the category of local systems if and
only if
$ r(e) = 0 $.

\medskip\noindent
{\bf 3.2.~Proposition} ({Non Lefschetz pencil case}).
{\it With the notation and the assumptions of {\rm (2.1)},
assume
$ Y = Y_{0} \mtims S $,
$ f = \pr_{2} $,
$ S $ is an affine rational curve, and
$ \pr_{1}\scirc i : X \to Y_{0} $ is birational, where
$ \pr_{i} $ is the
$ i $-th projection.
Assume furthermore that
$ Y_{0} $ has a nontrivial differential form of the
highest order.
Then {\rm (2.1.2)} does not split in the category of local
systems.
Furthermore, if there is a direct sum decomposition of
variations of Hodge structures
$ \cF' = \cF'_{1} \oplus \cF'_{2} $ such that
$ \cF'_{1} $ is a constant variation, then the extension
class between
$ \cF'' $ and
$ \cF'_{2} $ does not vanish in the category of local
systems.
}

\medskip\noindent
{\it Proof.}
Consider the short exact sequence
$$
\Gr_{m}^{W}H^{m}(Y,X;\bQ)\buildrel\alpha\over\to
\Gr_{m}^{W}H^{m}(Y,\bQ)\buildrel\beta\over\to
\Gr_{m}^{W}H^{m}(X,\bQ),
\leqno(3.2.1)
$$
where
$ W $ is the weight filtration of mixed Hodge structure [4].
Since
$ S $ is a smooth affine rational curve,
$ H^{1}(S,\bQ) $ has weights
$ > 1 $, and we have by the K\"unneth decomposition
$$
\Gr_{m}^{W}H^{m}(Y,\bQ) = H^{m}(Y_{0},\bQ).
$$
Then the last morphism
$ \beta $ in (3.2.1) is identified with the restriction
morphism
$$
H^{m}(Y_{0},\bQ) \to \Gr_{m}^{W}H^{m}(X,\bQ)
$$
by
$ X \to Y_{0} $, and its kernel is annihilated by the
restriction morphism to a nonempty open subvariety of
$ Y_{0} $, because
$ X \to Y_{0} $ is birational.

Thus the kernel of
$ \beta $ has level
$ < m $ (where the level of a Hodge structure is the difference
between the maximal and minimal numbers
$ p $ such that the
$ p $-th graded piece of the Hodge filtration does not vanish,
see [4]).
This implies that
$ \Im\, \beta \ne 0 $, because
$ H^{m}(Y_{0},\bQ) $ has level
$ m $ by the hypothesis on the highest form.
So we have a nonzero element
$$
u \in \Gr_{m}^{W}H^{m}(Y,\bQ) =
H^{0}(S,R^{m}f_{*}\bQ_{Y})
$$
such that
$ \beta(u) \ne 0 $.
Using the semisimplicity of polarizable Hodge structures [4],
we may assume that
$ \gamma_{m}(u) = 0 $ in
$ H^{0}(S,R^{m}g_{*}\bQ_{X}) $ (or equivalently, in
$ H^{m}(X_{s},\bQ) $ for
$ s \in S) $, because
$ H^{m}(X_{s},\bQ) $ has level
$ < m $.
Thus
$ u $ belongs to
$ H^{0}(S,\cF'') $.

For the first morphism
$ \alpha $ in (3.2.1), consider the Leray spectral sequence
in the category of mixed Hodge structures
$$
E_{2}^{p,q} = H^{p}(S,R^{q}h_{!}\bQ_{U}) \Rightarrow
H^{p+q}(Y,X;\bQ).
\leqno(3.2.2)
$$
This degenerates at
$ E_{2} $, because
$ E_{2}^{p,q} = 0 $ unless
$ p = 0 $ or
$ 1 $.
So we get the surjection
$$
\Gr_{m}^{W}H^{m}(Y,X;\bQ) \to
\Gr_{m}^{W}H^{0}(S,\cF)\,\,\,
\text{(see (2.1) for
$ \cF) $.}
\leqno(3.2.3)
$$

Now assume that (2.1.2) splits in the category of local
systems.
Then there exists
$ v \in H^{0}(S,\cF) $ whose image in
$ H^{0}(S,\cF'') $ is
$ u $.
Here we may assume
$ v \in \Gr_{m}^{W}H^{0}(S,\cF) $, because
$ H^{0}(S,\cF'') $ is pure of weight
$ m $, and the image commutes with
$ \Gr_{m}^{W} $.
But this contradicts the nonvanishing of
$ \beta(u) $ using the surjectivity of (3.2.3).
So the first assertion follows.

For the last assertion, it is enough to show the vanishing
of the extension class between
$ \cF''$ and
$ \cF'_{1} $ as local systems.
But this follows from (3.1.2) and (3.1.3)
because the first cohomology of a constant
variation of Hodge structures of weight
$ -1 $ on a rational curve has weight
$ > 0 $ (if it is nonzero).
This completes the proof of Proposition (3.2).

\medskip\noindent
{\bf 3.3.~Proposition} ({Surface case}).
{\it With the notation and the assumptions of {\rm (2.1)},
assume
$ m = 1 $ {\rm (}i.e.
$ Y $ is a surface and
$ X $ is a curve{\rm )},
$ S $ is an affine rational curve, and the fiber
$ Y_{s} $ of
$ f $ is not a rational curve.
Assume further that
$ f : Y \to S $ can be extended to a proper smooth morphism
$ f' : Y' \to S' $ so that
$ Y $,
$ S $ are open subvarieties of
$ Y' $,
$ S' $ respectively, and the closure
$ X' $ of
$ X $ in
$ Y' $ is nonsingular,
but it is not smooth over
$ S' $.
Then {\rm (2.1.2)} does not split in the category of
local systems.
}

\medskip\noindent
{\it Proof.}
In this case we have
$$
\cF'_{s} = \tH^{0}(X_{s},\bQ),\quad
\cF''_{s} = H^{1}(Y_{s},\bQ).
$$
Then in the notation of (3.1.3) we have
$$
H^{0}(S,\cH) = \Hom(\cF''_{\bQ},\cF'_{\bQ}) = 0,
$$
considering the monodromy of
$ \cF'_{s} $ around
$ S' \setminus S $.
So by (3.1.2) it is enough to show that the corresponding normal
function in (3.1.1) is nontorsion, and the assertion
is local on
$ S $ in the classical topology.

Let
$ s_{0} \in S' \setminus S $, and
$ \Delta $ be an open disk around
$ s_{0} $ in
$ S' $ such that
$ \Delta \cap S = \Delta^{*}\,
(:= \Delta\setminus\{s_{0}\}) $.
Then the assumption on
$ S' $ implies that there exist continuously
$ \Lambda_{s} = \{x_{s}, y_{s}\} \subset X_{s} $ for
$ s \in \Delta^{*} $ such that the action of the monodromy
around
$ s_{0} $ on
$ \Lambda_{s} $ is nontrivial (replacing
$ s_{0} $ if necessary).
Locally on
$ \Delta^{*} $, the difference
$ [x_{s}] - [y_{s}] $ defines an element
$ u_{s} \in \tH^{0}(X_{s}) $ and also a point
$ \xi_{s} $ of the Jacobian
$ J(Y_{s}) $ of
$ Y_{s} $.
Note that
$ \xi_{s} $ corresponds to the pull-back of the dual of
the short exact sequence of mixed Hodge structures
$ (2.1.2)_{s} $ by
$ u_{s} $ (using [2]), where
$ (2.1.2)_{s} $ is the stalk at
$ s $ of the exact sequence (2.1.2).

We have
$ \xi_{s} \ne 0 $ in
$ J(Y_{s}) $, because
$ Y_{s} $ is nonrational.
If a nonzero multiple of
$ \xi_{s} $ vanishes locally on
$ \Delta^{*} $, it defines a locally
constant section of the division points of the family of
Jacobians over
$ \Delta $ (because
$ f $ is proper smooth over
$ \Delta $), and we get a contradiction by considering the
monodromy around
$ s_{0} $.
So the assertion follows.

\medskip\noindent
{\bf 3.4.~Complement on the assumption of (3.3).}
The last assumption of Proposition (3.3) is satisfied in
the case of a generic Lefschetz pencil as follows.
Let
$ Y $ be a smooth surface embedded in a projective
space
$ \cP := \bP_{\bC}^{r} $
$ (r > 2) $.
Let
$ X $ be a (locally closed) smooth curve on
$ Y $.
Then there is a hyperplane of
$ \cP $ which intersects
$ Y $ transversely, but is tangent to
$ X $.

Indeed, let
$ \cP^{*} $ denote the dual projective space of
$ \cP $ consisting of hyperplanes
$ H $ of
$ \cP $, and
$ D_{Y} $,
$ D_{X} $ be the set of hyperplanes tangent to
$ Y $,
$ X $ respectively.
By definition,
$ D_{Y} $ is the image of a
$ \bP^{r-3} $-bundles
$ P_{Y} $ over
$ Y $ (where
$ P_{Y,y} $ consists of the hyperplanes tangent to
$ Y $ at
$ y $), and similarly for
$ X $ with
$ \bP^{r-3} $ replaced by
$ \bP^{r-2} $.
Let
$ x \in X $, and assume
$ P_{X,x}\subset D_{Y} $.
Then there exist a pencil of hyperplanes
$ \{H_{t}\}_{t\in \bP^{1}} $ contained in
$ P_{X,x} $ and a smooth analytic curve
$ C $ locally defined on
$ Y $ together with a nonconstant holomorphic map
$ \rho : C \to \bP^{1} $ such that
$ H_{\rho(c)} $ is tangent to
$ Y $ (and hence to
$ C $) at any
$ c \in C $.
But this implies the constancy of
$ \rho $, which is a contradiction.

\bigskip\bigskip
\centerline{\bf 4. Degeneration and Nearby Cycles}

\bigskip\noindent
In this section, we calculate the weight filtration
of Steenbrink, and prove Theorem (0.4) using the
results in Sections 2 and 3.

\medskip\noindent
{\bf 4.1.~Family of hypersurfaces.}
With the notation of (1.1), assume condition (0.1) in
Introduction.
Let
$ V_{1}, V_{2}, V_{3} $ be vector subspaces of
$ \cA_{Z}^{\delta}, \cA^{d-\delta}, \cA_{Z}^{d} $
such that
$ V_{1}, V_{3} $ generate
$ \cI_Z\otimes\cO_{Y}(\delta),\ \cI_Z\otimes\cO_{Y}(d) $
respectively and
$ V_{2} $ gives an embedding of
$ Y $ into a projective space.
We assume
$ V_{3} $ is the image of
$ V'_{3}\otimes V''_{3} $ where
$ V'_{3} $ is a vector subspace of
$ \cA_{Z}^{\delta} $ generating
$ \cO_Y(\delta) $ outside
$ Z $ and
$ V''_{3} $ is a vector subspace of
$ \cA^{d-\delta} $ giving an embedding of
$ Y $ into a projective space.
In the proof of Theorem (0.4),
$ V_{1}, V'_{3} $ and
$ V_{2}, V''_{3} $ will be respectively the restrictions of
$ \cA_{Z'}^{\delta} $ and
$ \cA^{d-\delta} $ defined for
some smooth projective variety
$ Y' $ containing
$ Y $ where
$ Z' \subset Y' $ is a subvariety whose intersection with
$ Y $ is
$ Z $.
(Actually, it is possible that
$ Z $ is empty, but
$ Z' $ is not.)
These are necessary to carry out an inductive
argument in (4.7).

We will identify
$ V_{1}, V_{2}, V_{3} $ with the corresponding affine spaces.
Let
$ \tV_{1}, \tV_{2}, \tV_{3} $ be their intersections with
$ \cV_{Z}^{\delta}, \cV^{d-\delta}, \cV_{Z}^{d} $.
Put
$ X_{P,Q} = X_{P}\cap X_{Q} $,
$ X_{P,Q,R} = X_{P}\cap X_{Q}\cap X_{R} $ for
$ P \in \cA^{i} $, etc.
Let
$$
S_{0} = \{(P,Q,R)\in \tV_{1}\mtims\tV_{2}\mtims \tV_{3}:
X_{P,Q}, X_{Q,R}, X_{P,Q,R} \,\,\text{are SCI}\},
$$
where SCI means smooth complete intersection.
By [19] (or (1.2)),
$ S_{0} $ is non empty, and
$ X_{P,R} $ has at most isolated singularities,
see also [13].
Here we can replace
$ S_{0} $ with a non empty subvariety because of (4.2) below.
Let
$$
S' = \{(P,Q,R,t)\in S_{0}\mtims\bC^{*} :
X_{PQ+tR} \,\, \text{is smooth}\},
$$
and
$ S'' = \cV_{Z}^{d}/\bC^{*} $.
We have a local system
$ L $ on
$ S'' $ whose stalks are given by
$ H^{m}(X,\bQ)_{\bot Z}^{\van} $ in Conjecture (0.2).
Consider a morphism
$ \rho : S' \to S'' $ which associates
$ PQ+tR $ to
$ (P,Q,R,t) $.
We apply the reduction argument in (4.2) below to this
so that the proof of Theorem (0.4) is reduced to the pull-back
of the local system to
$ S' $.
Here it is enough to show the {\it indecomposability} of
$ L $, because
$ L $ is semisimple by Deligne [4].

\medskip\noindent
{\bf 4.2.~Reduction argument.}
Let
$ L $ be a local system on a connected analytic space
$ S'' $.
Then it is simple (resp. indecomposable) if there exists a
morphism
$ \rho : S' \to S'' $ such that
$ \rho^{*}L $ is simple (resp. indecomposable).
This follows from the fact that the functor
$ \rho^{*} $ is exact and faithful (or using the corresponding
representation of the fundamental group).
A similar assertion also holds for the nearby cycle functor
$ \psi $.

\medskip\noindent
{\bf 4.3.~Semistable degeneration.}
With the notation of (4.1), there is a smooth projective morphism
$ f' : \cX' \to S' $ whose fiber
$ \cX_{s} := f^{\prime -1}(s) $ is
$ X_{PQ+tR} $ for
$ s = (P,Q,R,t) \in S' $.
Let
$ L' $ be a subsheaf of
$ R^{m}f'_{*}\bQ_{\cX'} $ whose stalk is the
orthogonal complement of the subspace generated by
$ H^{m}(Y,\bQ) $ together with the cycle classes of the
irreducible components of the closed subvariety
$ Z $.
Then
$ L' =\rho^{*}L $.
Let
$$
S = \{(P,Q,R,t)\in S_{0}\mtims\bC :
X_{PQ+tR}\,\, \text{is smooth with} \,\,
t \ne 0 \,\, \text{or} \,\, t = 0\},
$$
i.e.
$ S $ is the disjoint union of
$ S' $ and
$ S_{0} $.
Then
$ f' $ is naturally extended to
$ f : \cX \to S $.
However,
$ \cX $ has certain singularities.

Let
$ \cC $ be the closed subvariety of
$ \cX $ whose fiber over
$ (P,Q,R,t) \in S $ is
$ X_{Q,R} $.
Note that its restriction over
$ S' $ is a locally principal divisor.
Let
$ \pi : \tcX \to \cX $ be the blow-up along
$ \cC $.
Let
$ \tf : \tcX \to S $ denote the
composition with
$ f $.
We also denote by
$ t $ the function defined by the last component of
$ (P,Q,R,t) $.
Let
$ \tcX_{0}, \cX_{0}, S_{0} $ be the subvarieties of
$ \tcX, \cX, S $ defined by
$ t = 0 $ (this is compatible with the previous definition of
$ S_{0} $).
We have the induced morphism
$ \tf_{0} : \tcX_{0} \to S_{0} $.

Let
$ \cX_{P} $ be the closed subvariety of
$ \cX $ whose fiber over
$ (P,Q,R) \in S_{0} $ is
$ X_{P} $ (using the projection
$ \cX \to S_{0} $), and similarly for
$ \cX_{P,Q} $,
$ \cX_{P,Q,R} $, etc.
(Note that
$ \cX_{P} $,
$ \cX_{Q} $ are contained in
$ \cX_{0} $.)
Let
$ \cX_{P}^{\sim} $ be the blow-up of
$ \cX_{P} $ along
$ \cX_{P,Q,R} $.
Then
$ \tcX_{0} $ is a divisor with normal crossings on a smooth
variety
$ \tcX $, and its irreducible components are
$ \cX_{P}^{\sim} $ and
$ \cX_{Q} $, see [10].

Let
$ \psi $ denotes the nearby cycle functor, see [5].
Then
$$
R^{m}\tf_{*}\psi_{t}\bQ_{\tcX} =
\psi_{t}R^{m}\tf_{*}\bQ_{\tcX} \,
(= \psi_{t}R^{m}f'_{*}\bQ_{\cX'}),
\leqno(4.3.1)
$$
because the nearby cycle functor commutes with the direct image
under a proper morphism.

Let
$ W $ be the weight filtration on
$ \psi_{t}\bQ_{\tcX} $.
By Steenbrink [24], we have
$$
\aligned
\Gr_{m-1}^{W}\psi_{t}\bQ_{\tcX}
&= \bQ_{\cX_{P,Q}}[-1],
\\
\Gr_{m}^{W}\psi_{t}\bQ_{\tcX}
&= \bQ_{\cX_{Q}} \oplus
\bQ_{\cX_{P}^{\sim}},
\\
\Gr_{m+1}^{W}\psi_{t}\bQ_{\tcX}
&= \bQ_{\cX_{P,Q}}(-1)[-1],
\endaligned
\leqno(4.3.2)
$$
and
$ \Gr_{k}^{W}\psi_{t}\bQ_{\tcX} = 0 $ for
$ |k-m| > 1 $.

\medskip\noindent
{\bf 4.4.~Weight spectral sequence.}
Consider the weight spectral sequence in the category of local
systems
$$
{E}_{1}^{-k,j+k} = R^{j}\tf_{*}\Gr_{k}^{W}
\psi_{t}\bQ_{\tcX} \Rightarrow
\psi_{t}R^{j}f'_{*}\bQ_{\cX'}.
\leqno(4.4.1)
$$
This degenerates at
$ E_{2} $ by [24].
(In this case, it also follows from an easy calculation.)
The local system
$ \{H^{j}(X_{P})\} $ on
$ S_{0} $ can be identified with
$ R^{j}\tf_{*}\bQ_{\cX_{P}} $,
and similarly for
$ \{H^{j}(X_{P}^{\sim})\} $, etc., where the cohomology
groups are with rational coefficients unless otherwise
stated explicitly.
We have
$$
H^{j}(X_{P}^{\sim}) = H^{j}(X_{P}) \oplus H^{j-2}
(X_{P,Q,R})(-1),
$$
because
$$
\boR\pi_{*}\bQ_{\cX_{P}^{\sim}} =
\bQ_{\cX_{P}} \oplus \bQ_{\cX_{P,Q,R}}(-1)[-2],
\leqno(4.4.2)
$$
by the decomposition theorem [1].
Let
$ \cZ_{P}^{\sim} = \pi^{-1}(X_{P,Q,R}) $.
Then a canonical morphism
$ \bQ_{\cX_{P,Q,R}}(-1)[-2] \to \boR\pi_{*}\bQ_{\cX_{P}^{\sim}} $
and its right inverse (up to a sign) are given by the compositions
$$
\aligned
&\bQ_{\cX_{P,Q,R}}(-1)[-2] \to \boR\pi_{*}\bQ_{\cZ_{P}^{\sim}}
(-1)[-2] \to \boR\pi_{*}\bQ_{\cX_{P}^{\sim}},
\\
&\boR\pi_{*}\bQ_{\cX_{P}^{\sim}} \to
\boR\pi_{*}\bQ_{\cZ_{P}^{\sim}} \to
\bQ_{\cX_{P,Q,R}}(-1)[-2].
\endaligned
\leqno(4.4.3)
$$

The differential
$ d_{1} $ of the spectral sequence is induced by the
extension classes between the
$ \Gr_{k}^{W} $ which are given by the restriction and
Gysin morphisms.
Indeed, the extension class between
$ \Gr_{m-1}^{W} $ and
$ \Gr_{m}^{W} $ corresponds to a morphism
$ \Gr_{m}^{W} \to \Gr_{m-1}^{W}[1] $ in the derived
category of
$ \bQ $-modules on
$ \tcX $ (see [26]) and it is given by the restriction morphisms
$$
\bQ_{\cX_{Q}} \to \bQ_{\cX_{P,Q}},\quad
\bQ_{\cX_{P}^{\sim}} \to \bQ_{\cX_{P,Q}}.
$$
The direct image by
$ \pi $ of the last morphism is given by the restriction
and Gysin morphisms
$$
\bQ_{\cX_{P}} \to \bQ_{\cX_{P,Q}},\quad
\bQ_{\cX_{P,Q,R}}(-1)[-2] \to \bQ_{\cX_{P,Q}}
$$
(up to signs), using (4.4.2).
We have a similar assertion between
$ \Gr_{m}^{W} $ and
$ \Gr_{m+1}^{W} $.

Let
$ H^{m-1}(X_{P,Q})^{\van} $ denote the orthogonal complement
of the image of the (injective) restriction morphism
$ H^{m-1}(Y) \to H^{m-1}(X_{P,Q}) $ where
$ H^{m-1}(Y) $ can be replaced with
$ H^{m-1}(X_{P}) $ by the weak Lefschetz theorem,
and similarly for
$ H^{m-2}(X_{P,Q,R})^{\van} $ and also for
$ H^{m}(X_{P})^{\van} $.
Then we have the direct sum decompositions
$$
\aligned
&H^{m}(X_{P}) = H^{m}(X_{P})^{\van} \oplus
H^{m}(Y)^{\prim} \oplus H^{m-2}(Y)(-1),
\\
&H^{m}(X_{Q}) = H^{m}(X_{Q})^{\van} \oplus
H^{m}(Y)^{\prim} \oplus H^{m-2}(Y)(-1),
\\
&H^{m-2}(X_{P,Q,R}) = H^{m-2}(X_{P,Q,R})^{\van} \oplus
H^{m-2}(Y),
\endaligned
\leqno(4.4.4)
$$
together with the isomorphisms
$$
H^{m-2}(X_{P,Q})(-1) = H^{m}(X_{P,Q}) = H^{m-2}(Y)(-1).
$$
Here
$ H^{m}(Y)^{\prim} $ denotes the primitive cohomology, and
$ H^{m-2}(Y)(-1) $ in the first two decompositions of
(4.4.4) is actually the image of
$ H^{m-2}(Y)(-1) $ by the action of the cohomology class of
the hyperplane section.

We see that the stalk of
$ \Gr_{m}^{W}\psi_{t}R^{m}f'_{*}\bQ_{\cX'} $ is
the cohomology of the complex
$$
\aligned
H^{m-2}(X_{P,Q})(-1) \to H^{m}(X_{P}) \oplus H^{m}(X_{Q})
\oplus H^{m-2}(X_{P,Q,R})(-1)\,\,\,&
\\
\to H^{m}(X_{P,Q})&
\endaligned
\leqno(4.4.5)
$$
where the morphisms are induced by the restriction and Gysin
morphisms up to signs.
Using (4.4.4), the first morphism of (4.4.5) is
given by the identity on
$ H^{m-2}(Y)(-1) $ up to signs, and the last morphism
is given by the multiplications on
$ H^{m-2}(Y)(-1) $ by
$ d-\delta, \delta, d $ up to signs.
Indeed, the composition of
the Gysin and restriction morphisms
$$
H^{m-2}(X_{P,Q})(-1) \to H^{m}(X_{P}) \to H^{m}(X_{P,Q})
$$
coincides with the action of the restriction of the cycle class
of
$ X_{P,Q} $ in
$ X_{P} $.

Similarly we can verify
$$
\aligned
\Gr_{m-1}^{W}\psi_{t}R^{m}f'_{*}\bQ_{\cX'}
&= \{H^{m-1}(X_{P,Q})^{\van}\},
\\
\Gr_{m+1}^{W}\psi_{t}R^{m}f'_{*}\bQ_{\cX'}
&= \{H^{m-1}(X_{P,Q})^{\van}(-1)\}.
\endaligned
$$
Here the first assertion is easy, and the second follows
from it using duality.

Note that the extension classes between the
$ \Gr_{k}^{W} $ are induced by the restriction and
Gysin morphisms in the derived category of sheaves on
$ \tcX $ as above.

\medskip\noindent
{\bf 4.5.~Limit of the invariant part.}
Taking the nearby cycle functor (i.e. passing to the limit by
$ t \to 0 $), the restriction morphism
$ H^{m}(Y) \to H^{m}(X_{PQ+tR}) $ induces a morphism of
$ H^{m}(Y) $ to the middle term of the complex (4.4.5), which
is defined (up to signs) by using the isomorphism
$$
H^{m}(Y) = H^{m}(Y)^{\prim} \oplus H^{m-2}(Y)(-1)
\leqno(4.5.1)
$$
together with (4.4.4).
So
$ \Gr_{m}^{W} $ of the limit of
$ H^{m}(X_{PQ+tR})/H^{m}(Y) $ is given by
$$
(H^{m}(X_{P})^{\prim} \oplus
H^{m}(X_{Q})^{\prim} \oplus
H^{m-2}(X_{P,Q,R})^{\van}(-1))/
H^{m}(Y)^{\prim}.
\leqno(4.5.2)
$$
Note that
$$
H^{m}(X_{P})^{\prim} = H^{m-2}(X_{P})^{\van} \oplus
H^{m}(Y)^{\prim}
\leqno(4.5.3)
$$
(similarly for
$ H^{m}(X_{Q})^{\prim} $), and the division by
$ H^{m}(Y)^{\prim} $ is defined by using the diagonal
morphism.
In particular, the quotient (4.5.2) is isomorphic to
$$
H^{m}(X_{P})^{\prim} \oplus
H^{m}(X_{Q})^{\van} \oplus
H^{m-2}(X_{P,Q,R})^{\van}(-1),
\leqno(4.5.4)
$$
Here
$ \{H^{m}(X_{Q})^{\van}\} $ is an irreducible local
system on
$ S_{0} $.
Assuming an appropriate inductive hypothesis (see (4.7)
below), this also holds for the quotient of
$ H^{m-2}(X_{P,Q,R})^{\van} $ by the subspace generated
by the cycle classes of the intersection with
$ X_{Q} $ of the
$ (m/2) $-dimensional irreducible components of
$ Z $ if
$ m > 2 $.

\medskip\noindent
{\bf 4.6.~Cycle classes of the irreducible components.}
The cycle classes of the irreducible components of
$ Z $ in the limit of
$ H^{m}(X_{PQ+tR}) $ are given by using
$$
H_{\BM}^{m-2}(X_{P,R})(-1) \to H^{m}(X_{P}) \oplus
H^{m-2}(X_{P,Q,R})(-1),
\leqno(4.6.1)
$$
which is induced by the Gysin and restriction morphisms up to
signs.
(This can be verified by using (4.4.3).)
Here
$ H_{\BM}^{j}(V) $ denotes Borel-Moore
cohomology for an equidimensional variety
$ V $ in general, and is defined by
$ H^{j}(V,\bD_{V}(-n)[-2n]) $, where
$ \bD_{V} $ is the dualizing complex and
$ n = \dim V $.
If
$ V $ is compact, it is isomorphic to
$ H_{2n-j}(V)(-n) $.
The restriction of the dualizing complex to the smooth part
$ V_{\reg} $ is isomorphic to
$ \bQ_{V_{\reg}}(n)[2n] $, and we get the restriction
morphism
$ H_{\BM}^{j}(V) \to H^{j}(V') $ for any subvariety
$ V' $ of
$ V_{\reg} $.

The morphism to the second factor of (4.6.1) is
injective, i.e. the Gysin morphism
$ H^{m-2}(X_{P,Q,R})(-1)\to H^{m}(X_{P,R}) $ is surjective.
This follows from Artin's theorem (see [1]) which asserts
the vanishing of
$ H^{m}(X_{P,R}\setminus X_{P,Q,R}) $ because
$ X_{P,R}\setminus X_{P,Q,R} $ is affine and the constant
sheaf on it is semi-perverse up to a shift of complex by
$ m - 1 $.

We can verify that the limit of the cycle class of
$ Z $ is then defined by using the cycle class of the cycle in
$ X_{P} $ and the cycle class in
$ X_{P,Q,R} $ of the intersection of the
cycle with
$ X_{Q} $.

\medskip\noindent
{\bf 4.7.~Proof of Theorem (0.4).}
Assume first
$ m = 2 $ (and
$ Y' $ in (4.1) is
$ Y $).
In this case,
$ X_{P,Q,R} $ has dimension
$ 0 $ and is not connected.
We first fix
$ P, R $ and consider
$ \{X_{P,Q}\}_{Q} $ and
$ \{X_{P,Q,R}\}_{Q} $.
Let
$ Z_{i} \,(i>0) $ be the
$ 1 $-dimensional irreducible components of
$ Z $.
Then
$ X_{P,R} $ is the union of
$ \mcup_{i>0}Z_{i} $ and an irreducible curve
$ Z_{0} $, because
$ X_{P,R} \setminus Z $ is smooth and connected by (1.2)
(where
$ b = 0 $ if
$ \dim Z = 1 $), see also [13], [19].

We apply the dual of Proposition (2.7) or (3.3) to
$ \{X_{P,Q}\}_{Q} $ and
$ \{Z_{i,Q}\}_{Q} $ for
$ i \ge 0 $, where
$ Z_{i,Q} = Z_{i} \cap X_{Q} $.
Then we get the nontriviality of the extension class between
$ \Gr_{m-1}^{W} = \{H^{m-1}(X_{P,Q})^{\van}\} $ and
$ \{\tH^{m-2}(Z_{i,Q})(-1)\} $ for
$ i \ge 0 $, using (4.2).
Similarly we apply Proposition (2.7) or (3.2)
to show the nontriviality of the extension class between
$ \Gr_{m-1}^{W} $ and
$ \{H^{m}(X_{Q})^{\van}\} $, where we fix
$ Q, R $ or
$ P, R $ to apply the propositions.
We use Proposition (2.3) for the extension between
$ \Gr_{m-1}^{W} $ and any simple factor of
$ \{H^{m}(X_{P})^{\prim}\} $.
(Note that
$ \{H^{m}(X_{P})^{\prim}\} $ is semisimple by [4].)

Thus we get the nontriviality of the extension class
between
$ \Gr_{m-1}^{W}\psi_{t}L' $ and each simple factor of
$ \Gr_{m}^{W}\psi_{t}L' $ in the notation of (4.3), because
$ \Gr_{m}^{W}\psi_{t}L' $ is isomorphic to the quotient of
(4.5.4) by the image of the cycle classes of the
$ Z_{i}\,(i>0) $ (where the last term vanishes unless
$ \dim Z = 1 $), see (4.6).
We have the dual argument for the extension between
$ \Gr_{m}^{W} $ and
$ \Gr_{m+1}^{W} $.
Note that the
$ \{\tH^{m-2}(Z_{i,Q})\} $ are not isomorphic to each other
(considering the monodromy around
$ Q $ such that
$ X_{Q} $ is tangents to
$ Z_{j} $ at a smooth point for one
$ j $, but intersects
$ Z_{i} $ transversely at smooth points for
$ i \ne j $).
We also see that any simple factor of
$ \{H^{m}(X_{P})^{\prim}\} $ is not isomorphic to
$ \{\tH^{m-2}(Z_{i,Q})\} $, nor to
$ \{H^{m}(X_{Q})^{\van}\} $ (fixing
$ P $), and similarly between the last two.

Assume that there is a decomposition
$ \psi_{t}L' = L_{1}\oplus L_{2} $ in the notation of
(4.3).
We may assume that the
$ L_{i} $ are stable by the action of the monodromy
$ T $, because the decomposition is induced by that of
$ L' $.
In our case the weight filtration
$ W $ is defined by
$ W^{m-1} = \Im\,N $ and
$ W^{m} = \Ker\,N $.
Here
$ N = T -id $ because
$ N^{2} = 0 $.
If, for example,
$ L_{1} $ is not contained in
$ W^{m} $, then
$ NL_{1} = W^{m-1} $, because
$ NL_{1} \subset W^{m-1} $ is nonzero and
$ W^{m-1} = \Gr_{m-1}^{W} $ is simple.
In this case
$ L_{2} $ is contained in
$ W^{m} $, because otherwise it also contains
$ W^{m-1} $.
So the decomposition induces that of
$ \Gr_{m}^{W} $, and the nontriviality of the
above extension classes implies the triviality
of the decomposition.
Thus, using Deligne's semisimplicity theorem,
Theorem (0.4) is proved for
$ m = 2 $.

If
$ m = 1 $, then
$ X_{P} $ is smooth and connected by Theorem (1.2) with
$ r = 1 $, and
$ X_{P,Q,R} $ is empty.
So the argument is essentially same as above.

Assume now
$ m > 2 $.
We apply an inductive hypothesis to
$ X_{P,Q} $ and
$ Z \cap X_{Q} $ to show the irreducibility of
$ \{H^{m-2}(X_{P,Q,R})^{\van}\}_{R} $ divided by the image
of the cycle classes of the irreducible components of
$ Z $, see (4.6).
To carry out this induction we take general
$ P_{i},Q_{i},R_{i} $
such that
$ R_{i} = P_{i+1}Q_{i+1} + R_{i+1} $ and
$ Z $ is contained in
$ X_{P_{i}} $,
$ X_{R_{i}} $.
Note that the singular locus of
$ \mcap_{1\le j\le i}X_{P_j} $ has dimension
$ \le i - 2 $, see [19] (or (1.2)).
We apply the above argument to the restriction of
$ (P_{i}.Q_{i},R_{i}) $ to
$ Y'' := \mcap_{1\le j<i}X_{P_j,Q_j} $ for each
$ i $, and proceed by decreasing induction on
$ i $.
This is allowed by the definition of
$ S_{0} $ in (4.1) (where
$ Y', Y $ in (4.1) is
$ Y, Y'' $ here).
Here we use Proposition (2.7) or (3.2) (instead of (3.3))
to show the nontriviality of the extension between
$ \Gr_{W}^{m-1} $ and
$ \{H^{m-2}(X_{P,Q,R})^{\van}(-1)\} $ divided by the cycle
classes as above.
Note that even if
$ Y''\cap Z $ is empty, we still have some restriction to
$ (P_{i}.Q_{i},R_{i}) $ coming from
$ Z $.
At the first step of the induction, we have
$ m = 2 $ or
$ 1 $.
This completes the proof of Theorem (0.4).

\medskip\noindent
{\bf 4.8.~Generalization of Theorem (0.4).}
The assertion also holds for smooth zero loci of sections of
$ \cL_{1}\otimes\cL_{2} $ containing
$ Z $ where
$ \cL_{1}, \cL_{2} $ are line bundles such that
$ \cL_{1}\otimes\cI_{Z} $ is generated by its global sections,
and
$ \cL_{2} $ is very ample and satisfies one of the following
two conditions: either a general smooth hyperplane section of
$ \cL_{2} $ has a nontrivial differential form of the highest
degree, or the
$ 2 $-jets at each point is generated by the global sections of
$ \cL_{2} $ (e.g. it is the
$ (d -\delta) $-ple tensor of a very ample line bundle with
$ d \ge \delta+2 $).


\begin{thebibliography}{99}
\bigskip

\bibitem{1}
A. Beilinson, J. Bernstein and P. Deligne, Faisceaux pervers,
Ast\'erisque, vol. 100, Soc. Math. France, Paris, 1982.

\bibitem{2}
J. Carlson, Extensions of mixed Hodge structures, in Journ\'ees
de G\'eom\'etrie Alg\'ebrique d'Angers 1979, Sijthoff-Noordhoff
Alphen a/d Rijn, 1980, pp. 107--128.

\bibitem{3}
H. Clemens, Degeneration of K\"ahler manifolds, Duke Math. J.
44 (1977), 215--290.

\bibitem{4}
P. Deligne, Th\'eorie de Hodge I, Actes Congr\`es Intern. Math., 1970,
vol. 1, 425-430; II, Publ. Math. IHES, 40 (1971), 5--57; III ibid., 44
(1974), 5--77.

\bibitem{5}
\bysame, Le formalisme des cycles \'evanescents, in SGA7 XIII and XIV,
Lect. Notes in Math. vol. 340, Springer, Berlin, 1973, pp. 82--115 and
116--164.

\bibitem{6}
\bysame, La formule de Picard-Lefschetz, in SGA7 XV,
Lect. Notes in Math. vol. 340, Springer, Berlin, 1973, pp. 165-197.

\bibitem{7}
A. Dimca, Sheaves in Topology, Universitext, Springer, Berlin
(2003), to appear.

\bibitem{8}
A. Dimca and M. Saito, Monodromy at infinity and the weights of
cohomology, Compos. Math. 138 (2003), 55--71.

\bibitem{9}
D. Eisenbud, Commutative algebra with a view toward algebraic
geometry, Springer, New York, 1995.

\bibitem{10}
P. Griffiths and J. Harris, On the Noether-Lefschetz theorem and some
remarks on codimension two cycles, Math. Ann. 271 (1985) 31--51.

\bibitem{11}
L. Illusie, Autour du th\'eor\`eme de monodromie locale,
Ast\'erisque 223 (1994), 9--57.

\bibitem{12}
N. Katz, Etude cohomologique des pinceaux de Lefschetz, in Lect. Notes
in Math., vol. 340, Springer Berlin, 1973, pp. 254--327.

\bibitem{13}
S. Kleiman and A. Altman, Bertini theorems for hypersurface
sections containing a subscheme, Comm. Algebra 7 (1979), 775--790.

\bibitem{14}
S. Lefschetz, L'analysis situs et la g\'eom\'etrie alg\'ebrique,
Gauthier-Villars, Paris, 1924.

\bibitem{15}
J.D. Lewis, A survey of the Hodge conjecture (Second edition),
Monograph Series 10, American Mathematical Society, Providence,
RI, 1999.

\bibitem{16}
A.F. Lopez, Noether-Lefschetz theory and the Picard group of
projective surfaces, Mem. Amer. Math. Soc. 89 (1991).

\bibitem{17}
J. Milnor, Singular points of complex hypersurfaces, Ann. Math. Stud.
vol. 61, Princeton Univ. Press, 1968.

\bibitem{18}
A. Otwinowska,
Composantes de petite codimension du lieu de Noether-Lefschetz;
un argument en faveur de la conjecture de Hodge pour les
hypersurfaces, J. Algebraic Geom. 12 (2003), 307--320.

\bibitem{19}
\bysame, Monodromie d'une famille d'hypersurfaces, preprint.

\bibitem{20}
\bysame, Sur les vari\'et\'es de Hodge des hypersurfaces,
preprint (math.AG/0401092).

\bibitem{21}
M. Saito, Modules de Hodge polarisables, Publ. RIMS, Kyoto Univ., 24
(1988), 849--995.

\bibitem{22}
\bysame, Mixed Hodge Modules, Publ. RIMS, Kyoto Univ., 26
(1990), 221--333.

\bibitem{23}
\bysame, Admissible normal functions, J. Alg. Geom. 5 (1996),
235--276.

\bibitem{24}
J.H.M. Steenbrink, Limits of Hodge structures,
Inv. Math. 31 (1975/76), 229--257.

\bibitem{25}
J.H.M. Steenbrink and S. Zucker, Variation of mixed Hodge structure I,
Inv. Math. 80 (1985), 489--542.

\bibitem{26}
J.-L. Verdier, Cat\'egories d\'eriv\'ees, in SGA 4 1/2, Lect. Notes in
Math., Springer, Berlin, vol. 569, 1977, pp. 262--311.

\bibitem{27}
\bysame, Dualit\'e dans la cohomologie des espaces localement
compacts, S\'eminaire Bourbaki, No. 300.

\bibitem{28}
C. Voisin, Hodge theory and complex algebraic geometry, II,
Cambridge University Press, Cambridge, 2003.

\bibitem{29}
S. Zucker, Hodge theory with degenerating coefficients,
$ L_{2} $-cohomology in the Poincar\'e metric,
Ann. Math., 109 (1979), 415--476.

\end{thebibliography}
\end{document}